# FREQUENTIST OPTIMALITY OF BAYESIAN WAVELET SHRINKAGE RULES FOR GAUSSIAN AND NON-GAUSSIAN NOISE[1]

By Marianna Pensky

*University of Central Florida*

The present paper investigates theoretical performance of various Bayesian wavelet shrinkage rules in a nonparametric regression model with i.i.d. errors which are not necessarily normally distributed. The main purpose is comparison of various Bayesian models in terms of their frequentist asymptotic optimality in Sobolev and Besov spaces.

We establish a relationship between hyperparameters, verify that the majority of Bayesian models studied so far achieve theoretical optimality, state which Bayesian models cannot achieve optimal convergence rate and explain why it happens.

**1. Introduction.** Bayesian techniques for shrinking wavelet coefficients have become very popular in the last few years. Starting with the paper by Clyde, Parmigiani and Vidakovic [9], researchers turned to Bayesian methods in wavelet analysis (see, e.g., [1, 2, 3, 4, 6, 7, 8, 10, 17, 19, 24, 26]).

The majority of the above papers were devoted to the nonparametric wavelet regression model

$$Y_i = f(i/n) + Z_i, \qquad i = 1, \dots, n, \tag{1.1}$$

where $f$ is unknown and belongs to some functional space $\mathcal{F}$ and $Z_i$ are i.i.d. random variables. In most cases, the errors $Z_i$ were assumed to be normally distributed with zero mean and unknown variance $\sigma^2$. Some papers in engineering fields modeled errors as double-exponential (see, e.g., [20] or [23]); however, this approach was not very popular in statistical applications. The only paper known to the author which deals with nonnormal errors in a Bayesian framework is the one by Clyde and George [8] where the authors

Received February 2003; revised April 2005.
[1]Supported in part by NSF Grant DMS-05-05133.
*AMS 2000 subject classifications.* Primary 62G08; secondary 62C10.
*Key words and phrases.* Bayesian models, optimality, Sobolev and Besov spaces, nonparametric regression, wavelet shrinkage.

conduct comprehensive simulations with $Z_i$ having various distributions. In what follows we just assume that $Z_i$'s are i.i.d., symmetric and have finite fourth moment $EZ_1^4 < \infty$.

Now let us discuss the error model for wavelet coefficients which is used in Bayesian inference. If the $Z_i$'s in (1.1) are i.i.d. normal, the errors in wavelet coefficients are also i.i.d. normal, and most authors adopted this "real" model for Bayesian inference. However, some of the authors (Clyde and George [8] and Vidakovic [24]) were more adventurous and considered other distributions for error models, namely, distributions that result from a mixture of normal distributions over the scale parameter. Various choices of priors were suggested, and the posterior mean or median was used as an estimator of wavelet coefficients. The procedures were then justified by extensive simulations. Thus, a variety of algorithms have been suggested for applications; all of them have demonstrated good computational performance, but until the present time no attempts have been made to examine frequentist properties of these techniques or to compare various algorithms with each other.

To the best of the author's knowledge, three manuscripts have been produced recently with the objective to assess frequentist properties of Bayesian wavelet shrinkage and thresholding. The first one, by Abramovich, Amato and Angelini [1], explores optimality of Bayesian wavelet estimators in the case of normal errors in (1.1) and a combination of a point mass at zero and a normal density for the prior. The authors consider three choices of Bayesian wavelet estimators: the mean, the median and the mean thresholded on the basis of a Bayes factor rule [24], and then study optimality of the resulting wavelet regression estimators in Besov spaces. Johnstone and Silverman [19] go much further than Abramovich, Amato and Angelini [1], investigating theoretical and computational properties of the *empirical* Bayes thresholding rules. They also assume the errors to be normally distributed and choose priors to be a mixture of an atom probability at zero and a heavy-tailed density. The mixing weight, or sparsity parameter, for each level of the transform is chosen by maximizing marginal likelihood. The estimators can be based on a posterior mean or a posterior median, and they are optimal in Besov spaces. Finally, Autin, Picard and Rivoirard [5] compare a wide variety of Bayes shrinkage and thresholding rules based on both normal and heavy-tailed priors using a novel maxiset approach. However, although Autin, Picard and Rivoirard [5] and Johnstone and Silverman [19] develop breakthrough theory concerning regression with normally distributed errors, they do not discuss the case of a non-Gaussian distribution either for the actual error in (1.1) or as the error model in Bayesian inference.

The present paper is in a sense complementary to the three papers mentioned above as well as to Clyde and George [8]. Our interest was sparked by the excellent paper by Clyde and George [8], who present an amazing volume of simulations for various error distributions and various Bayesian models. All of their Bayesian models produce estimators as efficient (or better) than the ones based on traditional thresholding rules no matter whether the errors in (1.1) are normal or not. Yet, no theoretical study of these estimators is conducted; hence, some of the simulation results are stated but not explained.

The present paper examines various Bayesian models in terms of their frequentist asymptotic optimality without a predominant assumption that errors in (1.1) are normally distributed. Bayesian estimators for the wavelet coefficients are constructed using posterior means for a variety of error models and priors, and these models are compared in terms of their frequentist optimality. Although theoretical results of the paper are asymptotic, they allow one to explain some features of the simulations conducted by Clyde and George [8], namely, differences in the performances of "normal error–normal prior" model and the models with heavy-tailed error and prior distributions.

The material presented in the paper is interesting from several points of view. First, the relationship between hyperparameters is established, and it is verified that the majority of Bayesian models studied so far in the literature have indeed not only excellent computational properties but also theoretical optimality. Second, it is shown that some Bayesian models cannot achieve the optimal convergence rate and explained why this happens. Third, the prior distributions on wavelet coefficients are compared to see which of them agree with the assumption that the regression function belongs to a particular Besov space. The latter leads to a more comprehensive examination of the mixing weights. For example, we discover that the weights chosen in [1, 3] lead to convergence rates which differ from optimal (e.g., by the logarithmic factor in the case of normal errors and normal priors).

The rest of the paper is organized as follows. In Section 2 we introduce Bayesian models for wavelet coefficients. We use some "arbitrary" error model $\eta_j(\cdot)$ and a mixture of a point mass at zero and a density $\xi(\cdot)$ for a prior, keeping in mind that the actual distribution of wavelet coefficients is unknown at fine resolution levels and is asymptotically normal at coarse resolution levels according to the central limit theorem. In Section 3 we discuss assumptions on $\eta_j(\cdot)$ and $\xi(\cdot)$ and provide assertions about asymptotic optimality of regression estimators in Besov spaces for various choices of $\eta_j(\cdot)$ and $\xi(\cdot)$. In Section 4 we present comparison between various Bayesian models in tabular form and discuss theoretical results of Section 3. We also explain results of simulations conducted in [8]. Section 5 contains some auxiliary statements as well as proofs of the statements in Section 3.

**2. The Bayesian model.** Consider a standard nonparametric regression model (1.1) and assume that $f$ is square integrable and $Z_i$ are i.i.d. random variables with zero mean and finite variance. Application of the discrete wavelet transform (DWT) based on the wavelet $\psi$ to (1.1) yields

$$Y_{jk} = w_{jk} + \varepsilon_{jk}, \qquad j = 0, \dots, J-1, k = 0, \dots, 2^j - 1, J = \log_2 n, \tag{2.1}$$

where $\varepsilon_{jk}$ are uncorrelated random variables due to the unitary property of the DWT. We assume that the wavelet function $\psi$ has finite support and $s > r$ vanishing moments. If one is interested in the error of the resulting estimator of $f$ only at points $i/n$, one can just use periodized orthogonal wavelets. However, since in this paper we are interested in the $L^2$-norm based risk as the measure of the quality of the estimation procedure, some kind of boundary correction is necessary. For this reason, we shall use boundary coiflets introduced in [18, 19], which we describe in detail in Section 3.1.

Denote the scaling and the wavelet coefficients of the original function $f$ by $\tilde{\theta}_k$ and $\tilde{\theta}_{jk}$, respectively, so that $f$ can be reconstructed as

$$f(x) = \sum_{k=0}^{2^L-1} \tilde{\theta}_k 2^{L/2} \varphi(2^L x - k) + \sum_{j=L}^{\infty} \sum_{k=0}^{2^j-1} \tilde{\theta}_{jk} \psi_{jk}(x) \tag{2.2}$$

with $\tilde{\theta}_k = \int_{-\infty}^{\infty} 2^{L/2} \varphi(2^L x - k) f(x)\, dx$ and $\tilde{\theta}_{jk} = \int_{-\infty}^{\infty} \psi_{jk}(x) f(x)\, dx$. Here $\varphi(x)$ is a scaling function corresponding to the mother wavelet $\psi(x)$, $\psi_{jk}(x) = 2^{j/2}\psi(2^j x - k)$ and the value of $L$ will be defined later. Denote $\theta_{jk} = w_{jk}/\sqrt{n}$ and recall that $\tilde{\theta}_{jk} \approx \theta_{jk}$ (see, e.g., [25]). Later we shall provide a more detailed treatment of the relation between $\tilde{\theta}_{jk}$ and $\theta_{jk}$.

We shall use a Bayesian technique to construct estimators $\hat{\theta}_{jk}$ of $\theta_{jk}$ based on $Y_{jk}$. Since wavelet representations of a vast majority of functions contain only a few nonnegligible coefficients in their expansions, we place the following prior on the population discrete wavelet coefficient $w_{jk}$:

$$w_{jk} \sim \pi_{jn} \tau_{jn} \xi(\tau_{jn} x) + (1 - \pi_{jn}) \delta(0), \tag{2.3}$$

where $0 \le \pi_{jn} \le 1$ if $0 \le j \le J-1$ and $\pi_{jn} = 0$ if $j \ge J$, $\delta(0)$ is a point mass at zero and $w_{jk}$ are independent. The p.d.f. $\xi(x)$ is even and unimodal. Note that the majority of priors used previously for Bayesian wavelet inference follow the model (2.3). The factors $\pi_{jn}$ are the prior probabilities that a wavelet coefficient at level $j$ does not vanish. In what follows, however, we shall impose all conditions on the odds

$$\beta_{jn} = (1 - \pi_{jn})/\pi_{jn}. \tag{2.4}$$

Note that we allow dependence of $\pi_{jn}$ not only on the resolution level $j$ but also on $n$. It is most natural since the proportion of coefficients we are intending to keep depends not only on the function $f$ but also on the amount of data available: the larger $n$ is, the more reliable the estimators of coefficients are, the smaller the coefficients we can keep, the larger the value of $\pi_{jn}$ for a particular resolution level $j$.

Now let us discuss the distribution of errors. It follows from (1.1) that $\varepsilon_{jk} \approx n^{-1/2} 2^{j/2} \sum_{i=1}^{n} \psi(2^j i/n - k) Z_i$. Since $Z_i$ are i.i.d. with $EZ_1^4 < \infty$, the sequence $2^{j/2} n^{-1/2} \psi(2^j i/n - k) Z_i$ satisfies the Lyapunov condition with $\delta = 2$ provided $2^j/n \to 0$. Hence, if the resolution level is reasonably small ($j < J_0$ where $J - J_0 \to \infty$ as $n \to \infty$), the errors $\varepsilon_{jk}$ are asymptotically normally distributed $N(0, \sigma^2)$ and, thus, asymptotically independent, while at large resolution levels they are uncorrelated and have some unknown p.d.f. $\mu(x)$ which is symmetric and has a finite fourth moment. To incorporate both of these features into the distribution of $\varepsilon_{jk}$ we write it as

$$\varepsilon_{jk} \sim (1 - \lambda_j)(\sigma\sqrt{2\pi})^{-1} \exp(-x^2/(2\sigma^2)) + \lambda_j \mu(x). \tag{2.5}$$

Here $\sigma^2 = \int_{-\infty}^{\infty} x^2 \mu(x)\, dx = EZ_1^2$, $0 \le \lambda_j \le 1$, and $\lambda_j$ are close to 1 at large resolution levels and equal to zero at small resolution levels, namely, $\lambda_j = 1$ if $j = J - 1$ and $\lambda_j = 0$ if $j \le J_0$. For a more detailed treatment of asymptotic normality an interested reader can consult Neumann and von Sachs [21].

The difficulty of using (2.5) in Bayesian inference is that both $\mu(x)$ and $\lambda_j$ are unknown. For this reason, we shall choose the most general distribution for errors, namely,

$$\varepsilon_{jk} \sim \eta_j(x), \tag{2.6}$$

where $\eta_j$ are level-dependent symmetric densities. As we shall show later, one does not need the knowledge of the true distribution of $\varepsilon_{jk}$ and can obtain the optimal estimator of $f$ with a variety of error distributions $\eta_j$. Moreover, $\eta_j$ does not even need to be level dependent.

The model (2.6) generalizes the setup of Clyde and George [8], who considered error distributions other than normal in (1.1), namely, the distributions that can be obtained from the normal distribution by mixing over the scale parameter. It also follows the remarkable idea of Vidakovic [24], who assumed that the errors are normally distributed but used a double-exponential distribution for errors which resulted from the hierarchical Bayes approach, hence, conducting Bayesian inference with the error distribution different from the true one.

We shall conduct Bayesian inference for each wavelet coefficient separately. Denote

$$d_{jk} = Y_{jk}/\sqrt{n}, \qquad \nu_j = \sqrt{n}\tau_{jn}. \tag{2.7}$$

Taking into account the relation between $w_{jk}$ and $\theta_{jk}$ and (2.3)–(2.7), we derive that the posterior p.d.f. of $\theta_{jk}$ given $d_{jk}$ is of the form

$$p(\theta_{jk}|d_{jk}) = \frac{\sqrt{n}\eta_j(\sqrt{n}(\theta_{jk} - d_{jk}))\nu_j \xi(\nu_j \theta_{jk}) + \beta_{jn}\sqrt{n}\eta_j(\sqrt{n} d_{jk})\delta(0)}{\int_{-\infty}^{\infty} \sqrt{n}\eta_j(\sqrt{n}(x - d_{jk}))\nu_j \xi(\nu_j x)\, dx + \beta_{jn}\sqrt{n}\eta_j(\sqrt{n} d_{jk})}.$$

Choosing the posterior mean to be an estimator, we arrive at the following estimator of $\theta_{jk}$:

$$\hat{\theta}_{jk} = \frac{I_{1j}(d_{jk})}{I_{0j}(d_{jk}) + \beta_{jn}\sqrt{n}\eta_j(\sqrt{n}d_{jk})}, \qquad j \le J-1, \tag{2.8}$$

where

$$I_{ij}(d) = \int_{-\infty}^{\infty} x^i \sqrt{n}\eta_j[\sqrt{n}(x-d)]\nu_j \xi(\nu_j x)\, dx, \qquad i = 0, 1, \tag{2.9}$$

and $\hat{\theta}_{jk} = 0$ as $j \ge J$, so that the estimator of $f$ is of the form

$$\hat{f}(x) = \sum_{k=0}^{2^L-1} \hat{\theta}_k 2^{L/2}\varphi(2^L x - k) + \sum_{j=L}^{J-1}\sum_{k=0}^{2^j-1} \hat{\theta}_{jk}\psi_{jk}(x). \tag{2.10}$$

The objective of the present paper is to formulate conditions under which the estimator of $f$ based on coefficients $\hat{\theta}_{jk}$ is optimal and explore the cases where this optimality cannot be achieved. For any possible estimator $\tilde{f}$ of $f$ based on $n$ observations $X_1, \ldots, X_n$, define the mean integrated square error (MISE) over a function space $\mathcal{F}$ as

$$R_n(\mathcal{F}, \tilde{f}) = \sup_{f \in \mathcal{F}} E\|\tilde{f} - f\|^2_{L^2[0,1]}. \tag{2.11}$$

Donoho and Johnstone [14] showed that when the errors $Z_i$ in (1.1) are normally distributed and $f$ belongs to a ball $B^r_{p,q}(A)$ in the Besov space $B^r_{p,q}[0,1]$, then there exist constants $C_1$ and $C_2$ independent of $n$ such that

$$C_1 n^{-2r/(2r+1)} \le \inf_{\tilde{f}} R_n(B^r_{p,q}(A), \tilde{f}\,) \le C_2 n^{-2r/(2r+1)} \tag{2.12}$$

provided $r > \max(0, 1/p - 1/2)$ and $p, q \ge 1$. Since the Sobolev space $H^r = B^r_{2,2}$, the same rates of convergence hold for $H^r(A)$. Note that since the normal distribution for errors is a particular case of (2.5), the lower bounds in our situation cannot be smaller than (2.12). On the other hand, since for a majority of resolution levels ($j \le J_0$) the errors $\varepsilon_{jk}$ follow the normal model, we can expect to achieve convergence rates (2.12) as $n \to \infty$ for some choices of error models (2.6) and priors (2.3).

In what follows we shall compare various Bayesian models in terms of their ability to achieve optimal convergence rates (2.12) as $n \to \infty$.

## 3. Asymptotic optimality.

3.1. *Assumptions.* In what follows we shall formulate conditions on the p.d.f.'s $\xi(\cdot)$ and $\eta_j(\cdot)$ and parameters $\nu_j$, $\beta_{jn}$ and $J_0$. We warn the readers that not all of the conditions listed below will necessarily appear in every statement later.

Let $\phi$ and $\psi$ be boundary coiflets introduced in [18, 19], possessing $s > r$ vanishing moments and based on orthonormal coiflets supported in $[-S + 1, S]$, $s < S$. Assume that $p \geq 1$, $r > \max(1/2, 1/p)$ and

$$\begin{aligned} &r \geq r_p = 0.5[(1/p - 1/2) + \sqrt{(1/p-1/2)^2 + 2(1/p-1/2)}\,]I(p<2), \\ &L \geq \log_2(6S-6). \end{aligned} \tag{3.1}$$

Note that $r_p = 0$ when $p \geq 2$ and that $r_p \leq (1+\sqrt{5}\,)/4$ for any $p \geq 1$.

Let $\xi(x)$ and $\eta_j(x)$ be three times continuously differentiable everywhere except possibly $x = 0$, have finite fourth moments and satisfy the conditions

(3.2) (A1) $$|\xi^{(k)}(x)/\xi(x)| \leq C_\xi(1+|x|^{\lambda_\xi})^k, \qquad k = 1,2,3, \lambda_\xi \geq 0,$$

(3.3) (A2) $$|\eta_j^{(k)}(x)/\eta_j(x)| \leq C_\eta(1+|x|^{\lambda_\eta})^k, \qquad k = 1,2,3, \lambda_\eta \geq 0,$$

(3.4) (A3) $$|\eta_j(x)/\xi(x)| \leq C_{\xi,\eta},$$

(3.5) (A4) $$x_1^{2+\delta}\eta_j(x_1) \leq x_2^{2+\delta}\eta_j(x_2), \qquad x_1 > x_2 > C_\delta > 0, \delta \geq 0.$$

Condition (A4) just means that the functions $|x|^{2+\delta}\eta_j(x)$ are nonincreasing for sufficiently large $x$. Note that the constants $\lambda_\eta$, $C_\eta$, $C_{\xi,\eta}$, $\delta$ and $C_\delta$ are assumed to be independent of $j$ which requires some kind of uniformity of the p.d.f.'s $\eta_j$. The consequence of these restrictions is that asymptotic expansions of integrals $I_{ij}$ defined in (2.9) are valid with absolute constants independent of $j$, so in what follows we shall suppress the index $j$ in $I_{ij}, i = 0, 1$, unless this leads to confusion.

When conditions (A1) and (A2) hold with $\lambda_\xi = 0$ and $\lambda_\eta = 0$, following Johnstone and Silverman [17, 19], we shall say that $\eta_j$ and $\xi$ are *heavy-tailed* probability densities. The most common examples of the latter are p.d.f.'s of the double-exponential or Student $t$ distributions. In this situation, asymptotic expansions

$$I_0(d) \sim \nu_j\xi(\nu_j d), \qquad |I_1(d)/I_0(d) - d| = O(\nu_j/n) \qquad \text{if } \nu_j/\sqrt{n} \to 0, \tag{3.6}$$

$$I_0(d) \sim \sqrt{n}\eta(\sqrt{n}d), \qquad |I_1(d)/I_0(d)| = O(\sqrt{n}/\nu_j^2) \qquad \text{if } \nu_j/\sqrt{n} \to \infty, \tag{3.7}$$

which are proved later in Lemma 2, are valid for any $d$ as long as the relation between $\nu_j$ and $n$ holds. If $\lambda_\xi$ and $\lambda_\eta$ are positive, then the expansions (3.6) and (3.7) can be used under some restrictions on $d$ only. There is one important case though, when we can obtain asymptotic expansions for the integrals for any value of $d$: if $\xi$ and $\eta_j$ are normal p.d.f.'s and the variances of $\eta_j$ are bounded from above and from below by a positive value, then

$$|I_1(d)/I_0(d) - d| = O(|d|\nu_j^2/n) \qquad \text{if } \nu_j/\sqrt{n} \to 0, \tag{3.8}$$

$$I_0(d) \sim \nu_j \xi(\nu_j d) \qquad \text{if } \nu_j/\sqrt{n} \to 0, \tag{3.9}$$

$$|I_1(d)/I_0(d)| = O(n|d|/\nu_j^2) \qquad \text{if } \nu_j/\sqrt{n} \to \infty. \tag{3.10}$$

We denote

$$j_0 = (2r+1)^{-1} \log_2 n, \tag{3.11}$$

and assume that the parameter $\nu_j$ is of the form

$$\nu_j = C_1 2^{mj} \qquad \text{where } m = r + 1/2 \text{ if } p \geq 2. \tag{3.12}$$

The expression for $\nu_j$ when $1 \leq p < 2$ will be presented later [see (3.28)].

Remark 1. Assumptions about $\nu_j$ can be translated into the ones on $\tau_{jn}$ using relation (2.7), namely, $\tau_{jn} = C_1 2^{mj}/\sqrt{n}$. This coincides with the choice of $\tau_{jn}$ for the normal–normal model in [1, 3].

Remark 2. The assumption that $\varphi(x)$ and $\psi(x)$ are boundary coiflets as well as conditions (3.1) are introduced for the sake of obtaining convergence rates for a $L^2$-norm-based risk function. All statements of the paper will be true for $L = 0$ and an arbitrary $s$-regular scaling function $\varphi(x)$ and wavelet $\psi(x)$ with $s > \max(r, r + 1/2 - 1/p)$ if one replaces (2.11) by

$$R_n(\mathcal{F}, \tilde{f}) = \sup_{f \in \mathcal{F}} n^{-1} \sum_{i=1}^{n} E[\tilde{f}(i/n) - f(i/n)]^2.$$

Remark 3. Conditions on the existence of the fourth moment are purely technical conditions that are used for derivation of asymptotic expansions of the integrals $I_{ij}(\cdot)$. These conditions can be dropped and replaced by requesting that the conclusions of Lemma 2 remain valid. These conclusions, however, have to be verified individually for each combination of $\xi(\cdot)$ and $\eta_j(\cdot)$. In what follows we shall consider the case when $\eta_j$ and $\xi(\cdot)$ are p.d.f.'s of Student $t$ and Cauchy distributions, respectively. Later, however, we shall see that it is somewhat beneficial to apply a distribution with faster descent at $\pm\infty$ than the Cauchy.

3.2. *Optimality in Sobolev spaces and Besov spaces with $p \geq 2$.* It is well known that the Sobolev space $H^r$ can be characterized in terms of wavelet coefficients as $\sum_{j=L}^{\infty} \sum_{k=0}^{2^j-1} \tilde{\theta}_{jk}^2 (1 + 2^{2jr}) < \infty$ (see, e.g., [11], Section 9.2). Therefore, $f$ belongs to a ball $H^r(A)$ in $H^r$ space with $r > 1/2$:

$$f \in H^r(A) \quad \Longleftrightarrow \quad \sum_{j=L}^{\infty} \sum_{k=0}^{2^j-1} \tilde{\theta}_{jk}^2 (1 + 2^{2jr}) \leq A, \qquad r > \tfrac{1}{2}. \tag{3.13}$$

For the Besov ball $B^r_{p,q}(A)$ with $r > 1/p$ there is a similar relationship (see, e.g., [16]):

$$f \in B^r_{p,q}(A) \quad \Longrightarrow \quad \sum_{k=0}^{2^j-1} \tilde{\theta}^2_{jk} \le \begin{cases} B_1 2^{-2jr}, & \text{if } p \ge 2, \\ B_1 2^{-2j(r+1/2-1/p)}, & \text{if } 1 \le p < 2, \end{cases} \tag{3.14}$$

for some $B_1 > 0$. The cases $p \ge 2$ and $1 \le p < 2$ apply to spatially homogeneous and nonhomogeneous functions, respectively, and, as the reader will see later, the performance of Bayesian models varies greatly in those two cases. We shall start with the spatially homogeneous case $p \ge 2$ and study optimality of Bayesian estimators in this case. Since $H^r = B^r_{2,2}$, this is sufficient to study the general case of $f \in B^r_{p,q}(A)$ with $p \ge 2$; the results for a Sobolev ball will immediately follow.

In this subsection, we assume that $r > 1/2$ and $p \ge 2$. We also assume that $\nu_j$ is of the form (3.12) and $J_0$ are such that

$$2^{J_0} \ge n^{1/(2r)}. \tag{3.15}$$

Condition (3.15) is quite realistic and agrees with the assumption $J - J_0 \to \infty$ as $n \to \infty$ above, provided $r > 1/2$. In practice, the normality assumption can be checked via level-by-level testing.

In order to make understanding of a variety of assumptions and assertions of this section easier, we refer the reader to Table 1 in Section 4.2. In effect, we consider three kinds of models: (1) models with strictly optimal convergence rates (which include normal $\xi$–normal $\eta_j$ or heavy-tailed $\xi$–heavy-tailed $\eta_j$); (2) models optimal up to a log-factor (which include heavy-tailed $\xi$–normal $\eta_j$); (3) suboptimal models (which include normal $\xi$–heavy-tailed $\eta_j$). Within the first set of models we also study what happens when $\xi$ has a faster descent at $\pm\infty$ than $\eta_j$, that is, condition (A3) is violated.

3.2.1. *Models with strictly optimal convergence rates*: *normal* $\xi$*–normal* $\eta_j$ *or heavy-tailed* $\xi$*–heavy-tailed* $\eta_j$. For those models convergence rates are given by

THEOREM 1. *Let conditions* (A1) *and* (A2) *be valid with* $\lambda_\xi = \lambda_\eta = 0$, *or, alternatively, conditions* (3.8)–(3.10) *hold for any value of* $d$. *Let also conditions* (A3) *and* (A4) *hold. If for some positive* $\varepsilon$

$$\beta_{jn} = O((\sqrt{n}/\nu_j)^\alpha) \qquad \text{with } \alpha \le (2r+1)^{-1}(2+\delta) - 1 - \varepsilon \text{ as } j \le j_0, \tag{3.16}$$

*then*

$$R_n(B^r_{p,q}(A), \hat{f}\,) = O(n^{-2r/(2r+1)}), \qquad n \to \infty, p \ge 2. \tag{3.17}$$

Corollary 1. *If* $\eta_j(\cdot)$ *and* $\xi(\cdot)$ *are the p.d.f.'s of the double-exponential or the Student* $t$ *distributions, or both* $\xi(\cdot)$ *and* $\eta_j(\cdot)$ *are normal p.d.f.'s, satisfying conditions* (A2) *and* (A3) [*and* (A4) *for the* $t$ *distribution*], *then for* $\beta_{jn}$ *given by* (3.16) *the risk* $R_n(B^r_{p,q}(A), \hat{f})$ *is of the form* (3.17).

Note that in both Theorem 1 and Corollary 1 we assume that $\xi(x)$ and $\eta_j(x)$ have finite fourth moments. This is, however, a sufficient, not a necessary condition. Namely, the following statement is valid.

Corollary 1*. *If* $\eta_j(x)$ *are the p.d.f.'s of Student* $t_{2m_j+1}$ *distributions with* $m_j$ *integers,* $M_1 \le m_j \le M_2$, *and* $\xi(x)$ *is the p.d.f. of the Cauchy distribution, then* $R_n(B^r_{p,q}(A), \hat{f})$ *is of the form* (3.17) *provided* $\beta_{jn}$ *is given by* (3.16).

Now suppose condition (A3) is violated and $\xi(x)$ has faster descent to zero than $\eta_j(x)$ as $|x| \to \infty$. Then, in order for convergence rates to still hold, we need to impose extra conditions on $\beta_{jn}$, namely, $\beta_{jn}$ should be small for $j \le j_0$ (i.e., a priori more coefficients should be kept at low resolution levels). Consider an alternative assumption to (A3):

$$\text{(A3}^*\text{)} \qquad \begin{aligned} &|\eta_j(x)/\xi(x)| \le U(|x|) \\ &\text{where } U(x_1) \ge U(x_2) \text{ for any } x_1 > x_2 \ge 1. \end{aligned} \tag{3.18}$$

Theorem 2. *Let the conditions of Theorem* 1 *hold with* (A3) *replaced by* (A3*). *Let also for some* $\varepsilon > 0$

$$\begin{aligned} \beta_{jn} = O\Biggl(\min\Biggl(&\left[\frac{\sqrt{n}}{\nu_j}\right]^{(2+\delta)/(2r+1)-1-\varepsilon}, \\ &\left[\frac{\sqrt{n}}{\nu_j}\right]^{1/(2r+1)+\delta-\varepsilon} \frac{2^{-j/2}}{U(2B_1 2^{j/2})}\Biggr)\Biggr), \qquad j \le j_0, \end{aligned} \tag{3.19}$$

*where* $B_1$ *is the constant appearing in* (3.14). *Then* $R_n(B^r_{p,q}(A), \hat{f})$ *is of the form* (3.17).

Note that the constant $B_1$ is usually unknown. If $U(x)$ is a homogeneous function of some order [which happens, e.g., when $\eta_j(x)$ and $\xi(x)$ are the p.d.f.'s of $t$ distributions], then the value of $B_1$ has no bearing on $\beta_{jn}$. Otherwise, one can replace $B_1$ by any function of $n$ which grows infinitely as $n \to \infty$, for example, $\ln n$.

It is easy to see that condition (3.19) is more restrictive than (3.16). Condition (3.19) means that we a priori intend to keep many more coefficients than we "kill."

3.2.2. *Models optimal up to a log-factor*: *heavy-tailed $\xi$–normal $\eta_j$.* Unfortunately, for some combinations of distributions assumption (3.10) holds only when $\sqrt{n}|\sqrt{n}d|^{\lambda_\eta}/\nu_j \to 0$. This happens if $\eta_j(x)$ are normal distribution p.d.f.'s and $\xi(x)$ is a heavy-tailed p.d.f., for example, the double-exponential or Student $t$. In this situation we can ensure somewhat weaker conditions on $|I_1(d)/I_0(d)|$.

LEMMA 1. *If $\xi(x)$ and $\eta_j(x)$ are even unimodal p.d.f.'s such that $|I_1(d)| < \infty$ for any $j$ and $d$, then $|I_1(d)/I_0(d)| \le |d|$.*

Note that Theorem 1 and Corollaries 1 and 1* are valid for any distributions of $Z_i$'s in (1.1). However, if (3.10) does not hold for all values of $d$ and the errors in (1.1) are not normally distributed, one needs additional assumptions to achieve near the optimal convergence rate. These can be two different kinds of assumptions: either $\mu(x)$ in (2.5) should decrease reasonably fast at $\pm\infty$ or we should "kill" nearly all coefficients at the highest resolution levels. In particular, we define $\varrho_0(n)$ to be the unique positive solution of the equation

$$\int_{\varrho_0(n)}^{\infty} x^2 \mu(x)\, dx = n^{-2r/(2r+1)}, \tag{3.20}$$

where $\mu(x)$ and $\sigma$ are as in (2.5), and denote

$$\varrho(n) = \max[\varrho_0(n), \sigma\sqrt{2}\sqrt{2r(2r+1)^{-1}\ln n + 0.5\ln\ln n}\,]. \tag{3.21}$$

THEOREM 3. *Let conditions* (A1)–(A4) *hold with $\lambda_\xi = 0$ and $\lambda_\eta > 0$ and let $\beta_{jn}$ satisfy* (3.16). *If $\mu(\cdot)$ is such that*

$$\lim_{n\to\infty} n^{-1/4r}(\varrho(n))^{\lambda_\eta} = 0, \tag{3.22}$$

*or $\beta_{jn}$ increases quickly as $j \ge J_0$:*

$$\beta_{jn}^{-1} = O(\eta_j(2\sqrt{C_\mu} n^{r/(2r+1)})2^{rj/(2r+1)}), \qquad j \ge J_0, \tag{3.23}$$

*where $C_\mu = \int_0^\infty x^4 \mu(x)\, dx$, then*

$$R_n(B_{p,q}^r(A), \hat{f}\,) = O(n^{-2r/(2r+1)}(\ln n)^{\lambda_\eta/(2r+1)}), \qquad n \to \infty. \tag{3.24}$$

COROLLARY 2. *If $\eta_j(\cdot)$ are normal p.d.f.'s with variances bounded from above and from below by common positive numbers and $\xi(x)$ is the p.d.f. of the double-exponential or $t$ distribution, then* (3.24) *is valid provided* (3.16) *and either* (3.22) *or* (3.23) *holds.*

REMARK 4. Note that if the errors are normally distributed, then $\lambda_\eta = 1$ and $\varrho(n) = O(\sqrt{\ln n})$, so that condition (3.22) is unnecessary. If the errors follow the mixture model (2.5), condition (3.22) is valid whenever $\mu(x)$ decreases at an exponential rate or at a rate $|x|^{-a}$ with $a > 3 + 8r^2\lambda_\eta/(2r+1)$ as $|x| \to \infty$. The latter requires relatively fast descent (e.g., $a > 17/3$ for $r = 1$ and $\lambda_\eta = 1$) which may not be true. To overcome this difficulty, one can choose large values for $\beta_{jn}$ as $j \ge J_0$, suggested by (3.23). This measure will suppress coefficients at higher resolution levels and relax their dependence on a particular form of $\mu$.

3.2.3. *Suboptimal models: normal $\xi$–heavy-tailed $\eta_j$.* By now, we have addressed all possible models except for the ones where the prior $\xi(x)$ is a normal p.d.f. and the error distributions $\eta_j(x)$ are heavy-tailed. In this situation, assumption (3.8) may be invalid: instead of (3.8), by Lemma 1, we have

$$|I_1(d)/I_0(d) - d| = O(|d|) \qquad \text{if } \nu_j|\nu_j d|^{\lambda_\xi}/\sqrt{n} \ge M. \tag{3.25}$$

This happens when, for example, the $\eta_j(x)$ are p.d.f.'s of the double-exponential or the Student $t$ distribution. In this situation, convergence rates (3.17) or (3.24) do not hold.

THEOREM 4. *Let conditions* (A1) *and* (A2) *be valid with* $\lambda_\xi > 0$ *and* $\lambda_\eta = 0$. *Let conditions* (A3*) *and* (A4) *hold and* $\beta_{jn}$ *satisfy* (3.19). *Then*

$$R_n(B^r_{p,q}(A), \hat f\,) = O(n^{-2r/(2r+1+\lambda_\xi)}), \qquad n \to \infty. \tag{3.26}$$

Theorem 4 gives an upper bound for the risk. However, in order to conclude that certain combinations of distributions are asymptotically inferior to the others, we are interested in the lower bound for the risk.

COROLLARY 3. *If $\eta_j(x) = \eta(x)$ are identical p.d.f.'s of the double-exponential or the Student $t$ distribution and $\xi(x)$ is a normal p.d.f., then for some positive $C$*

$$R_n(H^r(A), \hat f\,) \ge Cn^{-2r/(2r+2)}, \qquad n \to \infty. \tag{3.27}$$

Observe that the lower bound in (3.27) is identical to the upper bound in (3.26) (since $\lambda_\xi = 1$ for the normal p.d.f.) and both are asymptotically larger than the optimal convergence rate (3.17). This is due to the fact that the *bias* of the estimator $\hat f$ converges at a slower rate. The latter happens because the model with the above combination of $\xi(\cdot)$ and $\eta(\cdot)$ fails to adapt to sparsity, namely, to the situation when at lower resolution levels we have very few relatively large coefficients. This is not surprising since in Corollary 3 we have a combination of a rather flat error model and a sharp prior p.d.f. which may fail to capture the actual value of the estimated coefficient.

3.3. *Optimality in Besov spaces with* $1 \le p < 2$. It follows from the previous section that the models which satisfy Corollary 1 attain the optimal convergence rate in Sobolev spaces and spatially homogeneous Besov spaces $B^r_{p,q}$, $p \ge 2$, with minimal assumptions whether the errors are normally distributed or not. Now, the question is whether all those models also achieve the optimal convergence rate in spatially nonhomogeneous spaces $B^r_{p,q}$, $1 \le p < 2$. Before giving an answer, we need to introduce new values of parameters for this case. Denote

$$j_1 = [(r + 1/2 - 1/p)(2r + 1)]^{-1} r \log_2 n,$$
$$J_0 = 0.5[\log_2 n + j_1],$$

and let $j_0$ and $\nu_j$ be defined by (3.11) and (3.12), respectively, with

$$m = \begin{cases} m_1 = r + 1/2 - 0.5(1/p - 1/2), & j \le j_0, \\ m_2 = (r + 1/2) - (1/p - 1/2)(1 + 1/r), & j_0 < j < j_1, \\ m_3 = r + 1/2, & j \ge j_1. \end{cases} \tag{3.28}$$

Observe that $m = r + 1/2$ for all resolution levels whenever $p = 2$. Note also that the resolution level $j_1$ is chosen so that

$$\sum_{j=j_1}^{\infty} \sum_{k=0}^{2^j-1} \theta_{jk}^2 = O(n^{-2r/(2r+1)}). \tag{3.29}$$

The restriction $r > r_p$ ensures that $J - j_1 \to \infty$, so that one can choose resolution level $J_0$, $j_1 < J_0 < J - 1$, such that $J - J_0 \to \infty$. Assume that

$$\text{(A5)} \qquad \int_{-\infty}^{\infty} z^2 \exp(-z^2/2\sigma^2)[\eta_j(1 + |z|)]^{-2}\, dz \le C, \tag{3.30}$$

where $\sigma$ is the TRUE standard deviation of the error and $C$ is a constant independent of $j$. Again we have the same three types of models as we had in the case of $p \ge 2$. However, neither of the first two types of models delivers strict optimality when $1 \le p < 2$. The third type of model is suboptimal even when $p \ge 2$, so we do not consider this type of model here.

3.3.1. *Normal* $\xi$*–normal* $\eta_j$ *or heavy-tailed* $\xi$*–heavy-tailed* $\eta_j$ *models.* Here we consider the subset of models studied in Section 3.2.1, namely, the models where the tails of $\xi$ are at least as heavy as the tails of $\eta_j$.

THEOREM 5. *Let* $f \in B^r_{p,q}(A)$ *with* $1 \le p < 2$ *and* $r > \max(r_p, 1/p)$. *Let* $\xi$ *and* $\eta_j$ *satisfy conditions* (A1) *and* (A2) *with* $\lambda_\eta = \lambda_\xi = 0$ *or let both of them be normal p.d.f.'s. Assume that* (A3) *and* (A4) *hold and*

$$\eta_j(x) \le C(x^2 + 1)^{-\nu/2} \exp(-\lambda |x|^\gamma). \tag{3.31}$$

*If* $\beta_{jn} = (\sqrt{n}/\nu_j)^\alpha$ *for* $j < j_1$ *where*

$$\alpha = \begin{cases} \alpha_1 \in (-\infty, \min(\nu - 3, (2r+2-1/p)^{-1}\nu - 1)), & \text{if } j \le j_0, \\ \alpha_2 \in \left[0, \nu - 2 + \dfrac{r+1-2/p}{2r(r+1)+1-p^{-1}(2r+2)}\right], & \text{if } j_0 < j \le j_1, \end{cases} \tag{3.32}$$

*whenever* $\gamma = 0$ *in* (3.31), *and* $\alpha \ge 0$ *whenever* $j_0 < j < j_1$ *and* $\gamma > 0$ *in* (3.31), *then*

$$R_n(B^r_{p,q}(A), \hat f\,) = O(n^{-2r/(2r+1)+\varepsilon_1}(\ln n)^{\varepsilon_2}), \qquad n \to \infty. \tag{3.33}$$

*Here*

$$\begin{array}{lll} \varepsilon_1 = 0, \qquad \varepsilon_2 = 8r/[\gamma(4r+p(2r+1))], & & \text{if } \gamma > 0, \\ \varepsilon_1 = [\alpha^*(1/p - 1/2)]/[\nu(2r+1)] > 0, & \varepsilon_2 = -p, & \text{if } \gamma = 0, \end{array} \tag{3.34}$$

*with* $\alpha^* = \max(1+\alpha_1, (1+\alpha_2)(2r+2)/r)$. *If* $\alpha_2 = 0$, *then* $\alpha^* = (2r+2)/r$, *the minimum possible value, and* $\varepsilon_1 = (2r+2)/[r\nu(2r+1)]$ *for* $\gamma = 0$.

Theorem 5 shows that the models where the $\eta_j$'s have exponential descent achieve optimal convergence rates up to a logarithmic factor while $\eta_j$ with polynomial descent lead to suboptimal convergence rates. Observe that condition (A5) does not exclude the normal model; it just requires that the variances of $\eta_j$ be larger than the actual variance of our data. In particular, if $\eta_j(x) = (\sqrt{2\pi\sigma_j^2})^{-1}\exp(-x^2/2\sigma_j^2)$, condition (A5) implies that $\sigma_j^2 > 2\sigma^2$. This does not contradict Abramovich, Amato and Angelini [1], who considered only the case of $\sigma_j \equiv \sigma$.

3.3.2. *Heavy-tailed $\xi$–normal $\eta_j$ models.* Now we study the models previously discussed in Section 3.2.2. These models under the assumption that $\eta_j$ are identical and the actual errors are normally distributed were investigated in depth by Johnstone and Silverman [19], who demonstrated strict optimality of the models in an empirical Bayes setup in a wide variety of Besov spaces. However, since Johnstone and Silverman [19] considered empirical Bayes estimators, we cannot just reproduce their results here.

THEOREM 6. *Let* $f \in B^r_{p,q}(A)$ *with* $1 \le p < 2$ *and* $r > \max(r_p, 1/p)$ *and let conditions* (A1)–(A4) *be valid. Let* $\eta_j(\cdot)$ *be normal p.d.f.'s and* $\lambda_\xi = 0$. *Assume also that either condition* (3.22) *or* (3.23) *is valid. If* $\beta_{jn} = (\sqrt{n}/\nu_j)^\alpha$ *where* $\alpha \ge 0$ *whenever* $j_0 < j < j_1$, *then*

$$R_n(B^r_{p,q}(A), \hat f\,) = O(n^{-2r/(2r+1)}(\ln n)^{\varepsilon_3}), \qquad n \to \infty, \tag{3.35}$$

*where* $\varepsilon_3 = \max(1/[2r+1], 4r/[4r+p(2r+1)])$.

3.4. *Optimality of the models with mixture distribution for errors.* In Sections 3.2 and 3.3 we considered only error distributions $\eta_j$ which satisfy certain assumptions uniformly. There is, however, one more interesting case, namely, that of $\eta_j(\cdot)$ being a mixture distribution mimicking (2.5):

$$\eta_j(x) = (1-\lambda_j^*)(\sqrt{2\pi}\sigma_0)^{-1}\exp(-x^2/(2\sigma_0^2)) + \lambda_j^*\zeta(x), \tag{3.36}$$

where $\lambda_j^* = 0$ for $j \le J_0$. Let $\zeta(x)$ be a heavy-tailed distribution (i.e., $\lambda_\zeta = 0$).

THEOREM 7. *Let $\lambda_\zeta = 0$ and $\xi(x) \ge C_{\xi,0}\exp(-x^2/(2\sigma_0^2))$ for any $x$. Let $\beta_{jn}$ satisfy the assumptions of Theorem 1 ($p \ge 2$) or Theorems 5 and 6 ($1 \le p < 2$). If $\xi(x)$ is a p.d.f. of the normal distribution, then as $n \to \infty$*

$$\begin{aligned} &R_n(B_{p,q}^r(A), \hat f\,) \\ &\quad = \begin{cases} O(n^{-2r/(2r+1)}), & \text{if } p \ge 2, \\ O(n^{-2r/(2r+1)}(\ln n)^{4r/(4r+p(2r+1))}), & \text{if } 1 \le p < 2, \end{cases} \end{aligned} \tag{3.37}$$

*the second assertion being valid provided $\sigma_0^2 > 2\sigma^2$, where $\sigma^2$ is the true error variance given by (2.5). If $\xi(x)$ is a heavy-tailed p.d.f. ($\lambda_\xi = 0$) and $\beta_{jn}$ and $\lambda_j^*$ are such that*

$$\beta_{jn} \ge \beta_0 \quad \text{and} \quad \lambda_j^* \le C_\lambda^*, \qquad j \ge J_0, \tag{3.38}$$

*then as $n \to \infty$*

$$R_n(B_{p,q}^r(A), \hat f\,) = O(n^{-2r/(2r+1)}(\ln n)^{\varepsilon}) \tag{3.39}$$

*where $\varepsilon = 1/(2r+1)$ if $p \ge 2$, and $\varepsilon = \varepsilon_3$ given by Theorem 6 if $1 \le p < 2$.*

The model (3.36) behaves exactly as the model with normal errors for $j \le J_0$. The advantage of adding a heavy-tailed term for $j \ge J_0$ is that even in the case when $\xi(x)$ is a heavy-tailed prior, no additional assumptions on the distribution $\mu(x)$ or $\beta_{jn}$ are necessary.

3.5. *Does $f \in B_{pq}^r$ a priori? More about prior odds $\beta_{jn}$.* So far, the choice of the error model $\eta_j(\cdot)$ was in the limelight. The main assertion about $\xi(\cdot)$ was that it should not have faster descent at $\pm\infty$ than $\eta_j(\cdot)$. However, it is $\xi(\cdot)$ that determines whether the regression function $f(x)$ belongs to a Besov space $B_{pq}^r$ a priori. Namely, the following statement is valid.

THEOREM 8. *If $\xi(\theta)$, $\nu_j$ and $\beta_{jn}$ are such that*

$$\int_{-\infty}^{\infty} |\theta|^{\max(p,q)}\xi(\theta)\,d\theta < \infty, \tag{3.40}$$

$$\lim_{n\to\infty}\sum_{j=L}^{\infty}[2^{j(r+1/2)}\beta_{jn}^{-1/p}\nu_j^{-1}]^{\min(p,q)} < \infty, \tag{3.41}$$

*then $f \in B_{p,q}^r$ a priori with probability 1.*

It is easy to see that condition (3.40) requires $\xi(\cdot)$ to have at least $\max(p,q) \geq 1$ finite moments, which immediately eliminates the Cauchy prior from the list. On the other hand, any prior with exponential descent ensures validity of (3.40).

COROLLARY 4. *Let $p \geq 2$, $q \geq 1$ and $\beta_{jn} = (\sqrt{n}/\nu_j)^{\alpha}$ with $\alpha > 0$ as $0 \leq j \leq j_0$. Let also*

$$\lim_{n\to\infty} \sum_{j=j_0+1}^{J-1} \beta_{jn}^{-\min(q/p,1)} < \infty. \tag{3.42}$$

*If the values of $\nu_j$ are given by (3.12) and $\xi(\theta)$ satisfies assumption (3.40), then $f \in B_{p,q}^r$ a priori with probability 1.*

Note that in the case when condition (A3) holds (Theorems 1 and 3), Corollary 4 can be applied. In order for one to be able to choose $\alpha > 0$ in Theorem 1, condition (A4) should hold with $\delta > 2r - 1$. Condition (3.42) is needed since Theorem 1 puts absolutely no restrictions on the values of $\beta_{jn}$ for $j \geq j_0$. The inequality (3.42) holds when, for example, $\beta_{jn} = (\nu_j/\sqrt{n})^{\alpha_1}$ with $\alpha_1 > 0$. The restriction (3.23) does not affect (3.42) since it imposes not small but large values of $\beta_{jn}$ for $j \geq J_0$. However, the case when condition (A3) is violated is troublesome since it calls for smaller values of $\beta_{jn}$ [see (3.19) in Theorem 2].

COROLLARY 5. *Let $1 \leq p < 2$ and $\nu_j$ and $\beta_{jn}$ be determined by Theorem 5 or 6. If (3.40) holds, $\beta_{jn} = (\sqrt{n}/\nu_j)^{\alpha}$ for $j < j_1$ where*

$$\alpha = \begin{cases} \alpha_1 \in (2p, \infty), & \text{if } j \leq j_0, \\ \alpha_2 \in (2p(r+1), \infty), & \text{if } j_0 < j \leq j_1, \end{cases} \tag{3.43}$$

*and*

$$\lim_{n\to\infty} \sum_{j=j_1+1}^{J-1} [\beta_{jn}^{-1/p}]^{\min(p,q)} < \infty, \tag{3.44}$$

*then $f \in B_{p,q}^r$ a priori with probability 1.*

Note that Corollary 5 is applicable under the conditions of Theorems 5–7 whenever assumption (3.31) on $\eta_j(\cdot)$ holds with $\gamma > 0$ or $\nu$ large enough, so that conditions (3.32) and (3.43) on $\alpha$ can be satisfied simultaneously. Again, assumption (3.44) is necessary since Theorems 5–7 put almost no restrictions on $\beta_{jn}$ when $j > j_1$.

Corollaries 4 and 5 imply that in order to ensure $f \in B_{p,q}^r$ a priori almost surely, one needs fairly large values $\beta_{jn}$, for example, such that the sums $\sum \beta_{jn}^{-\min(q/p,1)}$ are uniformly bounded [see (3.42)]. This fact motivates the choice $\pi_{jn} = 2^{-bj}$ in [1, 3], which is equivalent to $\beta_{jn} \sim 2^{jb}$. In the present paper, however, we choose $\beta_{jn} = (\sqrt{n}/\nu_j)^{\alpha}$, so that $\beta_{jn} \sim 1$ at the "optimal" resolution level $j_0$. Is this a reasonable choice? The answer is "yes" if one wants to obtain the optimal convergence rate for a wide range of models.

THEOREM 9. *Let conditions* (A1)–(A4) *be valid. Assume that* $p \geq 2$, $\nu_j$ *are defined by* (3.12) *and the densities* $\eta_j(x) = \eta(x)$ *are identical. If* $\pi_{jn} = O(2^{-bj})$ *with* $b > 0$, *then for some positive* $C$

$$R_n(B^r_{p,q}(A), \hat{f}\,) \geq C u(n) n^{-2r/(2r+1)} \quad \text{with } u(n) = [\eta^{-1}(n^{-b/(2r+1)})]^{2r/(r+1/2)} \to \infty. \tag{3.45}$$

*In particular,* $u(n) \sim (\ln n)^{2r/(2r+1)}$ *if* $\eta(x)$ *is a normal p.d.f.,* $u(n) \sim (\ln n)^{4r/(2r+1)}$ *if* $\eta(x)$ *is a double-exponential p.d.f. and* $u(n) \sim n^{br[(r+1/2)^2(2+\delta)]^{-1}}$ *if* $\eta(x) \sim |x|^{-(2+\delta)}$ *as* $|x| \to \infty$.

## 4. Discussion.

4.1. *Discussion of the simulations.* The case of the present paper is an occasion where extensive finite sample simulations have been carried out a few years before theoretical properties of the estimators have been studied. Clyde and George [8] discuss wavelet regression with errors being a scale mixture of normal errors and the priors on wavelet coefficients of the forms

$$\theta_{jk}|\lambda^*_{jk}, \gamma_{jk} \sim N(0, \sigma^2 c_j \gamma_{jk}/\lambda^*_{jk}), \qquad \lambda^*_{jk} \sim h^*, \qquad \gamma_{jk} \sim \text{Bernoulli}(\pi_j).$$

They consider threshold shrinkage and multiple shrinkage estimators, both based on the posterior mean, with the difference that the first one is calculated conditional on $\gamma_{jk}$ and data, while the second one is conditional on data only. Similarly to Johnstone and Silverman [17, 19], the authors estimate hyperparameters by maximizing marginal likelihood. They show that their estimators are computationally competitive with standard classical thresholding methods and also perform well in the case of nonnormally distributed errors. It is easy to check that the multiple shrinkage estimators of Clyde and George [8] coincide with the ones proposed in the present paper. Simulation study shows (see Figure 2 of their paper) that multiple shrinkage estimators are superior to the threshold shrinkage ones; therefore, the authors use multiple shrinkage estimators in their later simulations.

Clyde and George [8] suggest several choices for the prior and error distributions (normal prior–normal errors, normal prior–Student $t_\nu$ errors, $t_\nu$ prior–$t_\nu$ errors, etc.), but the only cases that made it to the actual simulations are normal prior–normal errors (N), Cauchy prior–$t_5$ errors (C5) and uncorrelated but dependent $t_5$ errors and prior (T5). Note that all of the three models satisfy Corollary 1 and are asymptotically optimal in $B^r_{pq}$ with $p \geq 2$ without additional assumptions on the error. Hence, though no studies for the less favorable cases are presented in the paper by Clyde and George [8], they probably have been eliminated as the ones producing inferior results.

Clyde and George [8] run simulations for the by now traditional test functions "blocks," "bumps," "doppler" and "heavisine" proposed by Donoho and Johnstone [12] and report results for three types of actual error distributions: normal, $t_5$ errors in wavelet domain and $t_5$ errors in data domain (see Figures 3, 5 and 6 of their paper). Before discussing Clyde and George's [8] findings we want to draw the reader's attention to the fact that their enumeration of resolution levels is just the opposite from ours: their highest resolution level is the most coarse while ours is the finest one.

Clyde and George [8] state that empirical Bayes wavelet estimators based on the models N, C5 and T5 are in the majority of cases superior to the traditional thresholding rules. However, they have no tool to explain the discrepancy in the performances of these three models, and discrepancies in the performance of various Bayesian models in general. We now are ready to make this final step.

While performing simulations with normal errors the authors discovered that model N gives somewhat better precision than the T5 and C5 models. The results reverse when the error has the Student $t_5$ distribution in the wavelet domain: models T5 and C5 lead to smaller MSE than model N. Clyde and George [8] remark that "performance of all the estimators tends to be worse under heavy-tailed error distributions" and that "the performance of N worsens the most." However, this is not surprising in view of the findings of Section 3.2. Note that heavy-tailed error distributions lead to increase in variance, and model N is the most vulnerable to this increase. Recall that at high resolution levels the variances of Bayesian estimators of the wavelet coefficients are bounded by a constant times $n/\nu_j^4$ for the models with $\lambda_\eta = 0$, while for the model N they are proportional to $n/\nu_j^2$ times the variance of the wavelet coefficient. This is the reason why the growth of the variance decreases efficiency of model N the most. The same phenomenon (although somewhat milder due to the central limit theorem) carries over to the case when the errors have $t_5$ distributions in the data domain.

The above also explains why Bayesian models show the least discrepancy for the "bumps" test function. This function, as Figure 4 of [8] shows, has the smallest value of prior variance which translates into higher $\nu_j$ in our notation. This leads to a higher proportion of the bias component in the overall error, and, by far, the highest value of the overall error among the four test functions. Hence, the variance component has the lowest weight in the overall error for the "bumps" function, so that all three models, N, C5 and T5, show similar performance in this case.

4.2. *Discussion and summary of the models.* Table 1 summarizes comparison of the Bayesian models carried out in the previous section. We assume that $\eta_j(\cdot) \equiv \eta(\cdot)$ are identical and consider three choices for $\xi(\cdot)$ and $\eta(\cdot)$: the normal, the double-exponential and the $t$ distribution.

The choices of $\xi$ and $\eta$ are listed in the first row and in the first column, respectively. The table presents asymptotic expressions for $n^{2r/(2r+1)} R(n, B^r_{pq}(A))$, which we denote by $\Delta_1$ when $p \geq 2$ and $\Delta_2$ when $1 \leq p < 2$; hence, $\Delta_1$ and $\Delta_2$ show deviation from the "ideal" rate $O(n^{-2r/(2r+1)})$. We also introduce the common parameter $\varsigma = 4r + p(2r+1)$. The case of the mixture distribution for errors is not covered by the table since it mimics the normal $\eta$–normal $\xi$ or normal $\eta$–heavy-tailed $\xi$, depending on whether $\xi$ is normal or heavy-tailed, with the only difference that the restriction (3.22) or (3.23) on the unknown distribution $\mu$ or $\beta_{jn}$, respectively, is unnecessary. To be more specific, we distinguish the situations when we have a mixed model (2.5) or a normal model [$\lambda_j \equiv 0$ in (2.5)] for errors. If both models give the same results with no additional assumptions, we leave the cell of the table unmarked; otherwise, if assumption (3.22) or (3.23) is required, we mark the cell with $\diamondsuit$. The cells where the additional assumption (3.19) on $\beta_{jn}$, $j \leq j_0$, is required are marked with $\triangle$.

Table 1 provides a comprehensive comparison of the models. As one can see, for spatially homogeneous Besov spaces ($p \geq 2$), the models of choice are the models which provide the optimal convergence rate, N–N, DE–DE, DE–T and T–T, where N, DE and T stand for normal, double-exponential and $t$ distributions, respectively; the first letter gives the choice of $\eta$ and the second that of $\xi$. These models provide optimal convergence rates with very few restrictions [only (A3) and (A4)], and no additional assumptions are required if the errors in (1.1) are not normally distributed. The models N–DE and N–T provide optimal convergence rates up to a logarithmic factor; however, the main flaw of these models is that they are sensitive to slow descent of the error distribution $\mu$ at $\pm\infty$: if (3.23) is not guaranteed, additional assumptions (3.22) on $\beta_{jn}$ for $j > J_0$ are necessary. To avoid this feature of N–DE and N–T models one can use models with the mixture error distribution (3.36), where $\zeta(\cdot)$ is a heavy-tailed p.d.f. We shall denote these models by N$\zeta$–N, N$\zeta$–DE and N$\zeta$–T, depending on the choice of $\xi$. Note that models N$\zeta$–DE and N$\zeta$–T behave similarly to N–DE and N–T but do not require assumption (3.22) or (3.23).

The situation changes if one considers spatially nonhomogeneous Besov spaces ($1 \leq p < 2$). In these cases, the model T–T becomes suboptimal while N–N, N–DE, N–T, DE–DE and DE–T still provide optimality up to a logarithmic factor. The model N–N, however, requires the very restrictive assumption $\sigma_0^2 > 2\sigma^2$, which can cripple empirical Bayes inference on parameters of the model. The models N–DE and N–T still require additional assumptions on $\beta_{jn}$ for $j > J_0$ in this case. Hence, the models of choice in this

TABLE 1
*Comparison of various models*

| $\eta(x)$ \ $\xi(x)$ | **Normal** $\frac{1}{\sqrt{2\pi\sigma_1^2}}\exp\left(-\frac{x^2}{2\sigma_1^2}\right)$ | **Double-exponential** $\frac{1}{2\sigma_1}\exp\left(-\frac{\lvert x\rvert}{\sigma_1}\right)$ | ***t*-distribution** $\frac{\Gamma((\nu_1+1)/2)(\nu_1\pi)^{-1/2}}{\Gamma(\nu_1/2)(1+(x^2/\nu_1))^{(\nu_1+1)/2}}$ |
|---|---|---|---|
| **Normal** $\frac{1}{\sqrt{2\pi\sigma_0^2}}\exp\left(-\frac{x^2}{2\sigma_0^2}\right)$ | $\Delta_1 = O(1)$ if $\sigma_0 \le \sigma_1$ or (3.19) <br> $\Delta_2 = O((\ln n)^{4r/\varsigma})$ if $\sqrt{2}\sigma < \sigma_0 \le \sigma_1$ | $\Delta_1 = O((\ln n)^{1/(2r+1)})$ <br> $\Delta_2 = O((\ln n)^{\kappa})$ $\kappa = \max\left(\frac{4r}{\varsigma}, \frac{1}{2r+1}\right)$ <br> ♢ | $\Delta_1 = O((\ln n)^{1/(2r+1)})$ <br> $\Delta_2 = O((\ln n)^{\kappa})$ $\kappa = \max\left(\frac{4r}{\varsigma}, \frac{1}{2r+1}\right)$ <br> ♢ |
| **Double-exponential** $\frac{1}{2\sigma_0}\exp\left(-\frac{\lvert x\rvert}{\sigma_0}\right)$ | $\Delta_1 = O(n^{2r/((2r+1)(2r+2))})$ <br> △ | $\Delta_1 = O(1)$ if $\sigma_0 \le \sigma_1$ or (3.19) <br> $\Delta_2 = O((\ln n)^{8r/\varsigma})$ if $\sigma_0 \le \sigma_1$ | $\Delta_1 = O(1)$ <br> $\Delta_2 = O((\ln n)^{8r/\varsigma})$ |
| ***t*-distribution** $\frac{\Gamma((\nu_0+1)/2)(\nu_0\pi)^{-1/2}}{\Gamma(\nu_0/2)(1+(x^2/\nu_0))^{(\nu_0+1)/2}}$ | $\Delta_1 = O(n^{2r/((2r+1)(2r+2))})$ <br> △ | $\Delta_1 = O(1)$ <br> △ | $\Delta_1 = O(1)$ if $\nu_0 \ge \nu_1$ or (3.19) <br> $\Delta_2 = O(n^{\varpi}(\ln n)^{-p})$ if $\nu_0 \ge \nu_1$ <br> $\varpi = \frac{2r+2}{r(2r+1)(1+\nu_0)}$ |

$\Delta_i = n^{2r/(2r+1)} R_n(B^r_{pq}(A), \hat{f})$ with $i = 1$ if $p \ge 2$ and $i = 2$ if $1 \le p < 2$.
$\varsigma = 4r + p(2r+1)$, $\sigma^2$ is the true error variance.
♢ requires assumption (3.22) or (3.23) for non-Gaussian errors in (1.1).
△ requires assumption (3.19) on $\beta_{jn}$, $j \le j_0$.

situation are DE–DE, DE–T (although they have a slightly higher power of the log-factor), N$\zeta$–DE, N$\zeta$–T, or N–DE and N–T if one is sure that the error distribution has fairly fast decay.

The T–DE model provides optimal convergence rates in spatially homogeneous Besov spaces but only under condition (3.19) which, as one can see from Corollary 4, prevents $f \in B^r_{pq}$ a priori. Hence, this model is inferior to DE–DE, DE–T and T–T and, for this reason, has not been studied in the case of $p < 2$. Note that N–N, DE–DE and T–T require condition (3.19) if $\sigma_0 > \sigma_1$ or $\nu_0 < \nu_1$.

The least advantageous choices in terms of asymptotic convergence rates are DE–N and T–N. Not only do these models fail to provide optimal convergence rates but they also require assumption (3.19) which prevents $f \in B^r_{pq}$ a priori.

Section 3.5 allows us to supplement our comparison with some restrictions on $\xi(\cdot)$. In order that $f \in B^r_{pq}$, condition (3.40) requires p.d.f. $\xi$ to have at least $\max(p, q)$ finite moments. Hence, condition (3.40) prevents $f \in B^r_{p\infty}$ whenever $\xi$ is a p.d.f. of the $t$ distribution. In general, Theorem 8 calls for prior distributions with faster descent at $\pm\infty$. However, it should be noted that the considerations of Theorem 8 have no bearing on the asymptotic convergence rates of the estimator.

Finally, some recommendations for a choice of a model can be drawn. If the errors in (1.1) are normally distributed, one can use DE–DE, DE–T, N–DE and N–T models which ensure optimality (up to a logarithmic factor) in both spatially homogeneous and nonhomogeneous Besov spaces. However, if the normality assumption is violated and errors may have a heavy-tailed distribution, DE–DE, DE–T, N$\zeta$–DE and N$\zeta$–T will be safer choices. In addition, if one wants to be true to the Bayesian spirit and make sure that $f \in B^r_{pq}$ a priori for any possible $q$, then the models where $\xi$ decreases exponentially (i.e., DE–DE and N$\zeta$–DE) may seem more preferable.

## 5. Proofs.

Proof of Lemma 1. Since $I_0(d)$ is an even and $I_1(d)$ is an odd function of $d$, consider the case $d > 0$. Note that

$$
\begin{aligned}
0 \le I_1(d) &= dI_0(d) - \int_{-\infty}^{\infty} (d-x)\sqrt{n}\eta[\sqrt{n}(d-x)]\nu_j\xi(\nu_j x)\,dx \\
&= dI_0(d) - \int_{-\infty}^{\infty} z\sqrt{n}\eta(\sqrt{n}z)\nu_j\xi(\nu_j(d-z))\,dz \le dI_0(d)
\end{aligned}
\tag{5.1}
$$

since the integral in (5.1) is nonnegative. Thus, $I_1(d) \le dI_0(d)$. $\square$

Proof of Theorem 1 is based on the following lemmas.

LEMMA 2. *If $\nu_j|d|$ is bounded or $\nu_j|\nu_j d|^{\lambda_\xi}/\sqrt{n}\to 0$, then as $\nu_j/\sqrt{n}\to 0$*

$$I_0(d)=\nu_j\xi(\nu_j d)[1+O(n^{-1}\nu_j^2|\nu_j d|^{2\lambda_\xi})], \tag{5.2}$$

$$\begin{aligned}\frac{I_1(d)}{I_0(d)}&=d-\eta_2\frac{\nu_j}{n}\frac{\xi'(\nu_j d)}{\xi(\nu_j d)}\left[1+O\left(\frac{\nu_j^2|\nu_j d|^{2\lambda_\xi}}{n}\right)\right]\\&=d-O\left(\frac{\nu_j|\nu_j d|^{\lambda_\xi}}{n}\right).\end{aligned} \tag{5.3}$$

*If $\sqrt{n}|d|$ is bounded or $\sqrt{n}|\sqrt{n}d|^{\lambda_\eta}/\nu_j\to 0$, then as $\sqrt{n}/\nu_j\to 0$*

$$I_0(d)\sim\sqrt{n}\eta(\sqrt{n}d)[1+O(n\nu_j^{-2}|\sqrt{n}d|^{2\lambda_\eta})], \tag{5.4}$$

$$\frac{I_1(d)}{I_0(d)}\sim\xi_2\frac{\sqrt{n}}{\nu_j^2}\frac{\eta'(\sqrt{n}d)}{\eta(\sqrt{n}d)}\left[1+O\left(\frac{n}{\nu_j^2}|\sqrt{n}d|^{2\lambda_\eta}\right)\right]=O\left(\frac{\sqrt{n}}{\nu_j^2}|\sqrt{n}d|^{\lambda_\eta}\right). \tag{5.5}$$

*Here $\eta_2=\int_{-\infty}^{\infty}x^2\eta(x)\,dx$, $\xi_2=\int_{-\infty}^{\infty}x^2\xi(x)\,dx$.*

PROOF. We shall give the proofs for (5.4) and (5.5); the proofs of (5.2) and (5.3) are conducted in a similar manner. Change variables $y=\nu_j x$ in (2.9) and use Taylor series expansion:

$$\begin{aligned}\frac{\nu_j^i}{\sqrt{n}}I_i(d)&=\int_{-\infty}^{\infty}y^i\eta\left(\sqrt{n}d-\frac{\sqrt{n}}{\nu_j}y\right)\xi(y)\,dy\\&=\int_{-\infty}^{\infty}y^i\Big[\eta(\sqrt{n}d)\\&\qquad-\frac{\sqrt{n}}{\nu_j}y\eta'(\sqrt{n}d)+\frac{n}{2\nu_j^2}y^2\eta''(\sqrt{n}d)\\&\qquad\qquad-\frac{n\sqrt{n}}{6\nu_j^3}y^3\eta'''(\sqrt{n}d)+\cdots\Big]\xi(y)\,dy.\end{aligned}$$

Letting $i=0$ and $i=1$ in (5.6) and taking into account that $\xi(\cdot)$ is an even function, we obtain (5.4) and (5.5). □

LEMMA 3. *Let $\varrho(n)$ be defined by (3.20) and (3.21) and let $\sigma$ be as in (2.5). Then*

$$E(d_{jk}-\theta_{jk})^{2i}=O(n^{-i}),\qquad i=1,2; \tag{5.6}$$

$$P(\sqrt{n}|d_{jk}-\theta_{jk}|>a\sqrt{\ln n})=o(n^{-a^2/(2\sigma^2)}),\qquad j<J_0; \tag{5.7}$$

$$E[(d_{jk}-\theta_{jk})^2I(\sqrt{n}|d_{jk}-\theta_{jk}|>\varrho(n))]=O(n^{-4r/(2r+1)}). \tag{5.8}$$

PROOF. Validity of statements (5.6) and (5.7) follows directly from the fact that [cf. (2.5)]

$$d_{jk} - \theta_{jk} \sim (1-\lambda_j)\sqrt{n}(\sqrt{2\pi}\sigma)^{-1}\exp\{-nx^2/(2\sigma^2)\} + \lambda_j\sqrt{n}\mu(\sqrt{n}x). \tag{5.9}$$

□

LEMMA 4. *Let $\phi$ and $\psi$ be boundary coiflets introduced in [18, 19] possessing $s > r$ vanishing moments and based on orthonormal coiflets supported in $[-S+1, S]$, $s < S$. Assume that $p \ge 1$, $r > \max(1/2, 1/p, r_p)$, where $r_p$ is defined in (3.1) and $L \ge \log_2(6S-6)$. If $f \in B^r_{p,q}(A)$, then for some absolute positive constants $A_1$ and $A_2$*

$$\sum_{j=L}^{J-1}\sum_{k=0}^{2^j-1}(\tilde{\theta}_{jk} - \theta_{jk})^2 \le A_1 n^{-2r/(2r+1)}, \tag{5.10}$$

$$\sum_{k=0}^{2^j-1}\theta_{jk}^2 \le A_2 2^{-2j[r-(1/p-1/2)_+]}. \tag{5.11}$$

PROOF. The proof is based on Proposition 5 of [18], which under conditions of Lemma 4 can be written as

$$\sum_{k=0}^{2^j-1}|\tilde{\theta}_{jk} - \theta_{jk}|^p \le AC(r,p,\varphi,\psi)2^{-rJp+(1/p-1/2)jp}. \tag{5.12}$$

To prove (5.10), consider cases $p \ge 2$ and $1 \le p < 2$ separately. If $p \ge 2$, the Hölder inequality leads to $\sum_{k=0}^{2^j-1}(\tilde{\theta}_{jk} - \theta_{jk})^2 = O(2^{-2rJ+2(1/p-1/2)j}2^{j(1-2/p)}) = O(2^{-2rJ})$. Adding the terms together we obtain

$$\sum_{j=L}^{J-1}\sum_{k=0}^{2^j-1}(\tilde{\theta}_{jk} - \theta_{jk})^2 = O\left(\sum_{j=L}^{J-1}2^{-2rJ}\right) = O(n^{-1}\ln n) = O(n^{-2r/(2r+1)}).$$

If $1 \le p < 2$, then $\sum_{k=0}^{2^j-1}(\tilde{\theta}_{jk} - \theta_{jk})^2 \le (\sum_{k=0}^{2^j-1}|\tilde{\theta}_{jk} - \theta_{jk}|^p)^{2/p} = O(2^{-2rJ+2(1/p-1/2)j})$. Summing the terms, we arrive at $\sum_{j=L}^{J-1}\sum_{k=0}^{2^j-1}(\tilde{\theta}_{jk} - \theta_{jk})^2 = O(2^{-2J(r+1/2-1/p)}) = O(n^{-2r/(2r+1)})$ provided $r > r_p$.

To derive (5.11), note that $\tilde{\theta} \in B^r_{p\infty}(A)$, where with some abuse of notation we use $B^r_{p\infty}$ to denote Besov space of infinite sequences. Hence $\sum_{k=0}^{2^j-1}\tilde{\theta}^2_{jk} = O(2^{2j(r-(1/p-1/2)_+)})$, which in combination with (5.10) yields (5.11). □

Proof of Theorem 1. Since the wavelet basis is orthonormal,

$$R_n(B^r_{p,q}(A),\hat{f}) = \sum_{k=0}^{2^L-1} E(\hat{\theta}_k - \tilde{\theta}_k)^2 + \sum_{j=L}^{J-1}\sum_{k=0}^{2^j-1} E(\hat{\theta}_{jk} - \tilde{\theta}_{jk})^2 + \sum_{j=J}^{\infty}\sum_{k=0}^{2^j-1} \tilde{\theta}^2_{jk}. \tag{5.13}$$

Observe that the first term in (5.13) is bounded by $2\sum_{k=0}^{2^L-1}[\mathrm{Var}(\hat{\theta}_k) + (\theta_k - \tilde{\theta}_k)^2] = O(n^{-1})$, while the last term is bounded by $A2^{-2rJ} = O(n^{-2r})$ due to (5.11). By Lemma 4, the second term in (5.13) is dominated by

$$\sum_{j=L}^{J-1}\sum_{k=0}^{2^j-1} E(\hat{\theta}_{jk} - \tilde{\theta}_{jk})^2 \le 2\sum_{j=L}^{J-1}\sum_{k=0}^{2^j-1} E(\hat{\theta}_{jk} - \theta_{jk})^2 + A_2 n^{-2r}, \tag{5.14}$$

that is, the main contribution to $R_n(B^r_{p,q}(A),\hat{f})$ is made by the first term in (5.14). Therefore, we need to construct an asymptotic upper bound for $R = \sum_{j=L}^{J-1}\sum_{k=0}^{2^j-1} E(\hat{\theta}_{jk} - \theta_{jk})^2 = R_1 + R_2$ with

$$R_1 = \sum_{j=L}^{j_0}\sum_{k=0}^{2^j-1} E(\hat{\theta}_{jk} - \theta_{jk})^2, \qquad R_2 = \sum_{j=j_0+1}^{J-1}\sum_{k=0}^{2^j-1} E(\hat{\theta}_{jk} - \theta_{jk})^2. \tag{5.15}$$

Let us examine each of the terms in turn. Denote

$$A_{jn}(d) = \beta_{jn}\sqrt{n}\eta_j(\sqrt{n}d)/I_0(d) \tag{5.16}$$

and note that $R_1 \le 2(R_{11} + R_{12})$ where

$$R_{11} = \sum_{j=L}^{j_0}\sum_{k=0}^{2^j-1} E\left(\frac{I_1(d_{jk})}{I_0(d_{jk})} - \theta_{jk}\right)^2, \qquad R_{12} = \sum_{j=L}^{j_0}\sum_{k=0}^{2^j-1} E\left(\frac{I_1(d_{jk})/I_0(d_{jk})}{1 + A_{jn}(d_{jk})} - \frac{I_1(d_{jk})}{I_0(d_{jk})}\right)^2. \tag{5.17}$$

To establish an asymptotic upper bound for $R_{11}$ observe that by combination of Lemma 2 and (3.8), for $j \le j_0$

$$\begin{aligned} &E(I_1(d_{jk})/I_0(d_{jk}) - \theta_{jk})^2 \\ &\quad \le 2[E(I_1(d_{jk})/I_0(d_{jk}) - d_{jk})^2 + E(d_{jk} - \theta_{jk})^2] \\ &\quad = O(E(\nu_j^4 d_{jk}^2/n^2) + \nu_j^2/n^2 + \sigma^2/n) \\ &\quad = O(\nu_j^4\theta_{jk}^2/n^2 + \nu_j^2/n^2 + \sigma^2/n), \end{aligned} \tag{5.18}$$

so that by (3.14)

$$R_{11} = O\left(\sum_{j=L}^{j_0}[2^{-2jr}\nu_j^4/n^2 + 2^j\nu_j^2/n^2 + 2^j/n]\right) = O(n^{-2r/(2r+1)}). \tag{5.19}$$

In the case of $R_{12}$, note that $R_{12} = R_{121} + R_{122}$ where

$$(5.20)\quad \begin{aligned} R_{121} &= \sum_{j=L}^{j_0}\sum_{k=0}^{2^j-1} E\biggl(\frac{I_1(d_{jk})}{I_0(d_{jk})}\frac{A_{jn}(d_{jk})}{1+A_{jn}(d_{jk})} I(\nu_j|d_{jk}| > C_\delta)\biggr)^2, \\ R_{122} &= \sum_{j=L}^{j_0}\sum_{k=0}^{2^j-1} E\biggl(\frac{I_1(d_{jk})}{I_0(d_{jk})}\frac{A_{jn}(d_{jk})}{1+A_{jn}(d_{jk})} I(\nu_j|d_{jk}| \le C_\delta)\biggr)^2. \end{aligned}$$

For $R_{121}$, recall that since $\nu_j \le \sqrt{n}$ and by assumptions (A3) and (A4)

$$(5.21)\quad \frac{\eta_j^2(\sqrt{n}d_{jk})}{\xi^2(\nu_j d_{jk})} \le C_{\xi,\eta}^2\biggl[\frac{(\sqrt{n}d_{jk})^{2+\delta}\eta_j(\sqrt{n}d_{jk})}{(\nu_j d_{jk})^{2+\delta}\eta_j(\nu_j d_{jk})}\biggr]^2 \frac{\nu_j^{4+2\delta}}{n^{2+\delta}} \le C_{\xi,\eta}^2 \frac{\nu_j^{4+2\delta}}{n^{2+\delta}}.$$

By (3.16), (5.16) and since $\delta + (2r+1)^{-1} > (2r+1)^{-1}(2+\delta) - 1$, the latter implies that

$$(5.22)\quad \begin{aligned} R_{121} &= O\Biggl(\sum_{j=L}^{j_0}\sum_{k=0}^{2^j-1} \beta_{jn}^2 \nu_j^{2+2\delta} n^{-(1+\delta)}[E(d_{jk}-\theta_{jk})^2 + \theta_{jk}^2]\Biggr) \\ &= O\Biggl(n^{-2r/(2r+1)} \sum_{j=L}^{j_0} \beta_{jn}^2 \biggl(\frac{\nu_j^2}{n}\biggr)^{1/(2r+1)+\delta}\Biggr) = O(n^{-2r/(2r+1)}). \end{aligned}$$

For the sake of construction of an upper bound for $R_{122}$, note that since $\nu_j|d_{jk}| \le C_\delta$ we have $I_0(d_{jk}) \sim \nu_j \xi(\nu_j d_{jk}) \ge \nu_j \xi(C_{dd})$. Thus, by $x/(1+x) \asymp \min(x,1)$, we derive that $A_{jn}(d_{jk}) = O(\min(1, \beta_{jn}^2 n \nu_j^{-2} \eta_j^2(\sqrt{n}d_{jk}))$. Let $\alpha_{jn}$ be a unique positive solution of $\beta_{jn}^2 n \nu_j^{-2}\eta_j^2(\alpha_{jn}) = 1$, that is,

$$(5.23)\quad \alpha_{jn} = \eta_j^{-1}\biggl(\frac{\nu_j}{\sqrt{n}\beta_{jn}}\biggr).$$

Since $\eta_j(x)$ are decreasing as $x > 0$, $A_{jn}(d_{jk}) \le 1$ iff $|d_{jk}| \ge \alpha_{jn}/\sqrt{n}$. Hence, by condition (A4), $E[d_{jk}^2 \min(1, A_{jn}^2(d_{jk}))] = E[d_{jk}^2 I(\sqrt{n}|d_{jk}| \le \alpha_{jn})] + E[d_{jk}^2 \beta_{jn}^2 n \nu_j^{-2} \eta_j^2(\sqrt{n}d_{jk}) I(\sqrt{n}|d_{jk}| > \alpha_{jn})]$, so that

$$(5.24)\quad E[d_{jk}^2 \min(1, A_{jn}^2(d_{jk}))] \le \alpha_{jn}^2/n + \beta_{jn}^2 n \nu_j^{-2} \eta_j^2(\alpha_{jn}^2)\alpha_{jn}^2/n = 2\alpha_{jn}^2/n,$$

since $x\eta_j(x)$ are nonincreasing functions of $x$ for $x > 0$. Now, note that, as $\nu_j|x| \ge C_\delta$, by condition (A4) $|x|^{2+\delta}\eta_j(x) \le C_\delta^{2+\delta}\eta_j(C_\delta) \le C_\delta^{2+\delta}\eta_j(0) \le C_\delta^{2+\delta} C_{\xi\eta}\xi(0) \equiv C_\eta^*$. Therefore, $\eta_j(x) \le C_\eta^{(1)}|x|^{-(2+\delta)}$ for any $x$ and $\eta_j^{-1}(z) \le C_\eta^{(2)} z^{-1/(2+\delta)}$. Let $\beta_{jn} = (\sqrt{n}/\nu_j)^\alpha$. Then, from (5.16)

$$(5.25)\quad \alpha_{jn} = \eta_j^{-1}([\nu_j/\sqrt{n}]^{1+\alpha}) = O([\sqrt{n}/\nu_j]^{(1+\alpha)/(2+\delta)}).$$

Combining (5.20), (5.24), (5.25) and Lemma 1, we derive

$$
\begin{aligned}
R_{122} &= O\left(\sum_{j=L}^{j_0}\sum_{k=0}^{2^j-1} E[d_{jk}^2 \min(1, A_{jn}^2(d_{jk}))]\right) \\
&= O\left(\sum_{j=L}^{j_0}\sum_{k=0}^{2^j-1} n^{-1}[n/\nu_j^2]^{(1+\alpha)/(2+\delta)}\right) \\
&= O\left(\sum_{j=L}^{j_0}\frac{2^j}{n}\left[\frac{n}{\nu_j^2}\right]^{(1+\alpha)/(2+\delta)}\right) \\
&= n^{-2r/(2r+1)} O\left(\sum_{j=L}^{j_0}\left[\frac{2^{(2r+1)j}}{n}\right]^{[1/(2r+1)-(1+\alpha)/(2+\delta)]}\right) \\
&= O(n^{-2r/(2r+1)})
\end{aligned}
\tag{5.26}
$$

since $(1+\alpha)(2r+1) < 2+\delta$. Combination of (5.19), (5.20), (5.22) and (5.26) ensures that $R_1 = O(n^{-2r/(2r+1)})$.

Now, consider $R_2$ given by (5.15). Since $|\hat{\theta}_{jk}| \le |I_1(d_{jk})/I_0(d_{jk})|$, by combination of Lemma 2 and (3.10)

$$
\begin{aligned}
R_2 &= O\left(\sum_{j=j_0+1}^{J-1}\sum_{k=0}^{2^j-1}\left[E\left(\frac{I_1(d_{jk})}{I_0(d_{jk})}\right)^2 + \theta_{jk}^2\right]\right) \\
&= O\left(\sum_{j=j_0+1}^{J-1}\sum_{k=0}^{2^j-1}\left[\frac{n}{\nu_j^4} + \frac{n^2 E d_{jk}^2}{\nu_j^4} + \theta_{jk}^2\right]\right) \\
&= O\left(\sum_{j=j_0+1}^{J-1}[2^j n 2^{-j(4r+2)} + 2^{-2jr}]\right) \\
&= O(n^{-2r/(2r+1)}),
\end{aligned}
\tag{5.27}
$$

which completes the proof. □

Proof of Corollary 1. Proof follows directly from Theorem 1 by Lemma 2 whenever $\lambda_\xi = \lambda_\eta = 0$. If $\xi(x)$ and $\eta_j(x)$ are normal p.d.f.'s, validity of (3.8)–(3.10) can be verified by direct calculation of $I_0(d)$ and $I_1(d)$. □

Proof of Corollary 1*. Using properties of the Fourier transform and formulae (8.432.5), (8.468), (3.944.5) and (3.944.6) of [15], one can show by direct calculation that $|I_1(d)/I_0(d) - d| = O(|d|\nu_j/\sqrt{n}), I_0(d) \sim \nu_j \xi(\nu_j d)$, if $\nu_j/\sqrt{n} \to 0$ and $|I_1(d)/I_0(d)| = O(\sqrt{n}|d|/\nu_j)$ if $\nu_j/\sqrt{n} \to \infty$, so Theorem 1 remains valid. □

Proof of Theorem 2. It is easy to see that condition (A3) is used only for derivation of the $R_{121}$ term in (5.20). Since $I(\nu_j|d_{jk}| > C_\delta) \le I(\sqrt{n}|d_{jk} -$ $\theta_{jk}| > a\sqrt{\ln n}\,) + I(\nu_j|d_{jk}| > C_\delta)I(\nu_j|d_{jk}| \le \nu_j|\theta_{jk}| + a\nu_j/\sqrt{\ln n}\sqrt{n}\,)$, we have $R_{121} = R_{1211} + R_{1212}$. Here, by Lemma 3,

$$\begin{aligned}R_{1211} &= O\Biggl(\sum_{j=L}^{j_0}\sum_{k=0}^{2^j-1} E[d_{jk}^2 I(\sqrt{n}|d_{jk}-\theta_{jk}| > a\sqrt{\ln n}\,)]\Biggr)\\ &= O\Biggl(\sum_{j=L}^{j_0}\sum_{k=0}^{2^j-1}[E(d_{jk}-\theta_{jk})^2 + \theta_{jk}^2 P(\sqrt{n}|d_{jk}-\theta_{jk}| > a\sqrt{\ln n}\,)]\Biggr)\\ &= O(n^{-2r/(2r+1)}) + O\Biggl(\sum_{j=L}^{j_0} 2^{-2jr} n^{-a^2/2\sigma^2}\Biggr) = O(n^{-2r/(2r+1)}),\end{aligned}$$

provided $a^2 \ge 4r\sigma^2/(2r+1)$. For $R_{1212}$, comparing with (5.21) and (5.22), take into account that $\nu_j|\theta_{jk}| \le B_1 2^{j/2}$ and $\nu_j\sqrt{\ln n}/\sqrt{n} = o(2^{j/2})$. Using (3.19), obtain

$$\begin{aligned}R_{1212} &= O\Biggl(\sum_{j=L}^{j_0}\sum_{k=0}^{2^j-1}\Biggl[\frac{\beta_{jn}^2}{\nu_j^2}\biggl(\frac{\nu_j^2}{n}\biggr)^{1+\delta}\biggl(\nu_j|\theta_{jk}| + a\frac{\nu_j}{\sqrt{n}}\sqrt{\ln n}\biggr)^2\\ &\qquad\qquad \times U^2\biggl(\nu_j|\theta_{jk}| + a\frac{\nu_j}{\sqrt{n}}\sqrt{\ln n}\biggr)\Biggr]\Biggr)\\ &= O\Biggl(\sum_{j=L}^{j_0}\frac{\beta_{jn}^2}{2^{2rj}}\biggl(\frac{\nu_j^2}{n}\biggr)^{1+\delta} 2^j U^2(2B_1 2^{j/2})\Biggr) = O(n^{-2r/(2r+1)}).\end{aligned}$$

□

Proof of Theorem 3. Note that since condition (3.10) is no longer valid, we need to derive new upper bounds for $R_2$. Note that $|\hat\theta_{jk} - \theta_{jk}| \le |\theta_{jk}| + |\hat\theta_{jk}|$ and $\sum_{j=j_0+1}^{J_0}\sum_{k=0}^{2^j-1}\theta_{jk}^2 = O(n^{-2r/(2r+1)})$, so we shall be concerned with the $|\hat\theta_{jk}|$ term only. Partition $\sum_j\sum_k\hat\theta_{jk}^2$ into the sum over $j_0 \le j \le J_0$ and $J_0 + 1 \le j \le J-1$ and denote the respective sums by $R_3$ and $R_4$.

To analyze $R_3$ note that $|\hat\theta_{jk}| \le |I_1(d_{jk})/I_0(d_{jk})|$ and that

$$\begin{aligned}1 &\le I(\sqrt{n}|\theta_{jk}| > a\sqrt{\ln n}\,) + I(\sqrt{n}|d_{jk}-\theta_{jk}| > a\sqrt{\ln n}\,)\\ &\quad + I(\sqrt{n}|d_{jk}-\theta_{jk}| \le a\sqrt{\ln n}\,)I(\sqrt{n}|\theta_{jk}| \le a\sqrt{\ln n}\,).\end{aligned}\tag{5.28}$$

Hence $R_3 = R_{31} + R_{32} + R_{33}$, where by Lemmas 1 and 3

$$R_{31} = O\Biggl(\sum_{j=j_0+1}^{J_0}\sum_{k=0}^{2^j-1} E[d_{jk}^2 I(\sqrt{n}|\theta_{jk}| > a\sqrt{\ln n}\,)]\Biggr)$$

$$= O\Biggl(\sum_{j=j_0+1}^{J_0}\sum_{k=0}^{2^j-1}\biggl[\theta_{jk}^2 + E(d_{jk}-\theta_{jk})^2 I\biggl(\frac{1}{n} < \frac{\theta_{jk}^2}{a^2\ln n}\biggr)\biggr]\Biggr)\tag{5.29}$$

$$= o(n^{-2r/(2r+1)}),$$

$$R_{32} = O\left(\sum_{j=j_0+1}^{J_0}\sum_{k=0}^{2^j-1} E[d_{jk}^2 I(\sqrt{n}|d_{jk}-\theta_{jk}| > a\sqrt{\ln n}\,)]\right)$$

$$= O\left(\sum_{j=j_0+1}^{J_0}\sum_{k=0}^{2^j-1}[\sqrt{E(d_{jk}-\theta_{jk})^4}\right.$$

$$\left.\times\sqrt{P(|d_{jk}-\theta_{jk}|\sqrt{n} > a\sqrt{\ln n}\,)} + \theta_{jk}^2]\right) \tag{5.30}$$

$$= O\left(\sum_{j=j_0+1}^{J_0}\sum_{k=0}^{2^j-1}[n^{-1}n^{-2r/(2r+1)} + \theta_{jk}^2]\right)$$

$$= O(2^{J_0}n^{-(4r+1)/(2r+1)}) = o(n^{-2r/(2r+1)})$$

provided $a^2 \geq 8\sigma^2 r/(2r+1)$. To derive an asymptotic upper bound for $R_{33}$, note that

$$\begin{aligned} &I(\sqrt{n}|d_{jk}-\theta_{jk}| \leq a\sqrt{\ln n}\,)I(\sqrt{n}|\theta_{jk}| \leq a\sqrt{\ln n}\,)\\ &\quad\leq I(\sqrt{n}|d_{jk}| \leq 2a\sqrt{\ln n}\,)\\ &\quad\leq I(\sqrt{n}|d_{jk}| \leq 2a\sqrt{\ln n}\,)I\left(\frac{\sqrt{n}}{\nu_j}(\sqrt{\ln n}\,)^{\lambda_\eta} \to 0\right)\\ &\qquad+ I\left(\frac{\nu_j}{\sqrt{n}}(\sqrt{\ln n}\,)^{-\lambda_\eta} = O(1)\right)\\ &\quad\leq I(\nu_j^{-1}\sqrt{n}(\sqrt{n}|d_{jk}|)^{\lambda_\eta} \to 0) + I(2^{j(2r+1)} = O(n(\ln n)^{\lambda_\eta})). \end{aligned} \tag{5.31}$$

Note that Lemma 2 and $\sqrt{n}\nu_j^{-1}(\sqrt{n}|d_{jk}|)^{\lambda_\eta} \to 0$ imply that

$$\begin{aligned} &E[[I_1(d_{jk})/I_0(d_{jk})]^2 I(\nu_j^{-1}\sqrt{n}(\sqrt{n}|d_{jk}|)^{\lambda_\eta} \to 0)]\\ &\quad= O(n\nu_j^{-4}(\sqrt{n}|d_{jk}|)^{2\lambda_\eta} I(\nu_j^{-1}\sqrt{n}(\sqrt{n}|d_{jk}|)^{\lambda_\eta} \to 0))\\ &\quad= O(\nu_j^{-2}) = O(2^{-(2r+1)j}). \end{aligned}$$

Therefore, by calculations similar to (5.27), the portion of $R_{33}$ corresponding to the first term in (5.31) is $O(n^{-2r/(2r+1)})$.

By Lemma 1, $E[I_1(d_{jk})/I_0(d_{jk})]^2 = O(E[d_{jk}-\theta_{jk}]^2 + \theta_{jk}^2)$, so the second term in the portion is

$$O\left(\sum_{j=j_0+1}^{J_0}\sum_{k=0}^{2^j-1}[n^{-1}+\theta_{jk}^2]I[2^j = O(n^{1/(2r+1)}(\ln n)^{\lambda_\eta/(2r+1)})]\right)$$

$$= O(n^{-2r/(2r+1)}(\ln n)^{\lambda_\eta/(2r+1)}).$$

Consequently,

$$R_{33} = O(n^{-2r/(2r+1)}(\ln n)^{\lambda_\eta/(2r+1)}), \tag{5.32}$$

and (5.29), (5.30) and (5.32) imply that $R_{321} = O(n^{-2r/(2r+1)}(\ln n)^{\lambda_\eta/(2r+1)})$.

Derivation of an asymptotic upper bound for $R_4$ is, by and large, similar to that for $R_3$. First, assume that condition (3.22) holds. Rewrite (5.28) and (5.31) with $a\sqrt{\ln n}$ replaced by $\varrho(n)$ and observe that under assumption (3.22) the second indicator $I(2^{j(2r+1)} = O(n(\ln n)^{\lambda_\eta}))$ in (5.31) vanishes. Now, same as before, $R_{41} = o(n^{-2r/(2r+1)})$ and $R_{42} = O(n^{-2r/(2r+1)})$ by (5.8). To ensure that $R_{43} = O(n^{-2r/(2r+1)})$ repeat calculations for the $R_{33}$ corresponding to the first term in (5.31).

If assumption (3.22) is violated, the only term in the proof which is affected is $R_{43}$. Namely, the second indicator in (5.31) does not vanish. Recall that $R_{43} = O(\sum_{j=J_0+1}^{J-1}\sum_{k=0}^{2^j-1} E[\hat\theta_{jk}^2 I(\sqrt{n}|d_{jk}| \le 2\varrho(n))])$. Since $|\hat\theta_{jk}| \le |I_1(d_{jk})|/[\beta_{jn}\sqrt{n}\eta_j(\sqrt{n}d_{jk})]$ and $|I_1(d)| \le \eta(0)\sqrt{n}\nu_j^{-1}\int_{-\infty}^{\infty}|z|\xi(z)\,dz$, we derive that, by (3.23), $R_{43} = O(\sum_{j=J_0+1}^{J-1}\sum_{k=0}^{2^j-1} E[\beta_{jn}\nu_j\eta_j(2\varrho(n))]^{-2})$. To complete the proof note that the assumption that $\mu$ has at least four finite moments implies that $\varrho(n) \le \sqrt{C_\mu}n^{r/(2r+1)}$. Consequently, $R_{43} = O(n^{-2r/(2r+1)})$. □

Proof of Corollary 2. The proof follows directly from Lemma 2. □

Proof of Theorem 4. The fact that condition (3.25) replaces (3.8) affects the term $R_{11}$ only. Denote $\Delta_{jk} = I_1(d_{jk})/I_0(d_{jk}) - I_1(\theta_{jk})/I_0(\theta_{jk})$ and observe that

$$R_{11} \le 2\sum_{j=L}^{J-1}\sum_{k=0}^{2^j-1}[E(I_1(\theta_{jk})/I_0(\theta_{jk}) - \theta_{jk})^2 + E\Delta_{jk}^2] \equiv 2(R_{111} + R_{112}). \tag{5.33}$$

For an upper bound on $R_{111}$ note that (5.11) implies that $|\theta_{jk}| \le \sqrt{A}2^{-jr}$, so that $(\nu_j|\theta_{jk}|)^{\lambda_\xi}\nu_j/\sqrt{n} = O(2^{j(r+\lambda_\xi/2+1/2)}/\sqrt{n})$. The last expression does not turn to zero only if $2^j > Cn^{1/(2r+1+\lambda_\xi)}$ for some $C > 0$. Also, it follows from Lemma 2 that

$$(\nu_j|\theta_{jk}|)^{\lambda_\xi}\nu_j/\sqrt{n} \to 0 \quad\Longrightarrow\quad |I_1(\theta_{jk})/I_0(\theta_{jk}) - \theta_{jk}| = O(1/\sqrt{n}). \tag{5.34}$$

Hence,

$$\begin{aligned} R_{11} &= O\left(\sum_{j=L}^{j_0} 2^j/n\right) + O\left(\sum_{j=L}^{j_0}\sum_{k=0}^{2^j-1}\theta_{jk}^2 I[2^j > Cn^{1/(2r+1+\lambda_\xi)}]\right) \\ &= O(n^{-2r/(2r+\lambda_\xi+1)}). \end{aligned}$$

To examine $R_{112}$ consider three separate cases. If both $(\nu_j|\theta_{jk}|)^{\lambda_\xi}\nu_j/\sqrt{n}$ and $(\nu_j|d_{jk}|)^{\lambda_\xi}\nu_j/\sqrt{n}$ tend to zero, then $\Delta_{jk} \le |I_1(d_{jk})/I_0(d_{jk}) - d_{jk}| + |d_{jk} - \theta_{jk}| + |I_1(\theta_{jk})/I_0(\theta_{jk}) - \theta_{jk}|$ and, by Lemma 2 and (5.34) we have $E\Delta_{jk}^2 = O(1/n)$. Consequently, the respective portion of $R_{112}$ is $O(n^{-2r/(2r+1)})$. If $(\nu_j|\theta_{jk}|)^{\lambda_\xi}\nu_j/\sqrt{n}$ does not tend to zero, then, as was mentioned earlier, $2^j > Cn^{1/(2r+1+\lambda_\xi)}$ for some $C > 0$. Hence the respective portion of $R_{112}$ is $O(\sum_{j=L}^{j_0}\sum_{k=0}^{2^j-1}[\theta_{jk}^2 I(2^j > Cn^{1/(2r+1+\lambda_\xi)}) + E(d_{jk} - \theta_{jk})^2]) = O(n^{-2r/(2r+\lambda_\xi+1)})$. The third case occurs if $(\nu_j|\theta_{jk}|)^{\lambda_\xi}\nu_j/\sqrt{n} \to 0$ but $(\nu_j|d_{jk}|)^{\lambda_\xi}\nu_j/\sqrt{n}$ does not tend to zero. Then $\Delta_{jk}^2 = O(1/n) + O(|d_{jk} - \theta_{jk}|^2) + O(d_{jk}^2)$ and $\sqrt{n}|d_{jk} - \theta_{jk}| \ge C^*(n/\nu_j^2)^{(\lambda_\xi+1)/2\lambda_\xi}$. Considering cases $n/\nu_j^2 > \ln n$ and $n/\nu_j^2 < \ln n$ separately and applying Lemma 3, we obtain $R_{112} = O(n^{-2r/(2r+\lambda_\xi+1)})$. □

The proof of Corollary 5 is based on the following lemmas.

LEMMA 5. *Let $\xi(\cdot)$ and $\eta_j(\cdot) \equiv \eta(\cdot)$ be such that $|I_1(d)/I_0(d) - d| \ge |d|/2$ for $d = \sqrt{A/2}\,2^{-i_0 r}$ with $i_0 < j_0$. Then $R(n, H^r(A)) \ge (A/8)2^{-2i_0 r}$.*

PROOF. Consider $f(x) = \theta_{i_0,0}^*\psi_{i_0,0}(x)$ with $\theta_{i_0,0}^* = \sqrt{A/2}\,2^{-i_0 r}$, $i_0 \ge L$. It is easy to check that $\theta_{i_0,0} = \theta_{i_0,0}^*$ and $\theta_{jk} = 0$ otherwise. Thus, the coefficients $\theta_{jk}$ satisfy condition (5.11) with $A_2 = A$ and $p = 2$. For this function $f$, under the condition of Lemma 5, the bias term exceeds $R_1 \ge C_2^2\sum_{j=L}^{j_0}\sum_{k=0}^{2^j-1}\theta_{jk}^2 = \theta_{i_0,0}^2 \ge (A/8)2^{-2i_0 r}$. □

LEMMA 6. *If $\xi(x)$ is a p.d.f. of the normal distribution and $\eta(x)$ is a p.d.f. of the double-exponential or the Student $t$ distribution, then Lemma 5 holds with $i_0 = (2r+2)^{-1}\log_2 n$.*

PROOF. Let $\nu_j/\sqrt{n} \to 0$, and $d = \sqrt{A/2}\,2^{-jr}$. If $\eta(x)$ is the Student $t$ distribution, then direct calculations show that $|I_1(d)/I_0(d) - d| \sim |d(C_\nu\sigma^2/(\nu_j d)^2 - 1)|$ where $C_\nu$ is a constant depending on the degrees of freedom of the $t$ distribution; hence $|I_1(d)/I_0(d) - d| \ge |d|/2$ for $j \ge i_0$ with $2^{i_0} \ge n^{1/(2r+2)}$.

In the case where $\eta(x)$ is the double-exponential distribution, $|I_1(d)/I_0(d) - d| \sim |d - \sigma^2\sqrt{n}\nu_j^{-2}| \ge |d|/2$ for $j \ge i_0$ with $i_0$ satisfying $2^{i_0} > (8n\sigma^4A^{-1})^{1/(2r+1)}$. □

PROOF OF COROLLARY 3. Corollary 3 follows directly from Lemmas 5 and 6. □

The proof of Theorem 5 is based on the following lemma.

LEMMA 7. *If $f \in B^r_{p,q}(A)$, then for $b > 0$*

$$\sum_j \sum_{k=0}^{2^j-1} \theta_{jk}^2 I(\sqrt{n}|\theta_{jk}| \le (\ln n)^b) = O(n^{-2r/(2r+1)}[\ln n]^{4br/(2r+1)}),$$

$$\sum_j \sum_{k=0}^{2^j-1} n^{-1} I(\sqrt{n}|\theta_{jk}| > (\ln n)^b) = O(n^{-2r/(2r+1)}[\ln n]^{-pb}).$$

PROOF. The proof is an adaptation of the proof in Donoho, Johnstone, Kerkyacharian and Picard [13]. □

PROOF OF THEOREM 5. Similarly to (5.15), partition $R$ as $R = R_1 + R_5 + R_2$ where $R_1$, $R_5$ and $R_2$ correspond to $j \le j_0$, $j_0 < j \le j_1$ and $j > j_1$, respectively, and consider each term separately. Note also that in condition (A4), $\delta$ can be made as large as one desires if $\gamma > 0$ and $\delta \le \nu - 2$ if $\gamma = 0$.

*Low resolution levels*: $j \le j_0$. Note that $\nu_j$ is chosen so that $\nu_j/\sqrt{n} = o(1)$ and $\sum_{j=L}^{j_0} \sum_{k=0}^{2^j-1} (n^{-2}\nu_j^2 + n^{-2}\nu_j^4 \theta_{jk}^2) = O(n^{-2r/(2r+1)})$, so that the upper bound (5.19) for $R_{11}$ is still valid. Partition $R_{12} = R_{121} + R_{122}$ as in (5.20). Consider the cases $\gamma > 0$ and $\gamma = 0$ separately. If $\gamma > 0$, then taking into account that $\nu_j |d_{jk}| > C_\delta$,

$$\frac{\eta_j^2(\sqrt{n} d_{jk})}{\xi^2(\nu_j d_{jk})} \le \exp\left\{-\lambda\left[\left(\frac{n}{\nu_j^2}\right)^{\gamma/2} - 1\right] C_\delta^\gamma\right\} = O\left(\left[\frac{\nu_j}{\sqrt{n}}\right]^{\upsilon}\right) \tag{5.35}$$

for any $\upsilon > 0$. In the case $\gamma = 0$,

$$\frac{\eta_j^2(\sqrt{n} d_{jk})}{\xi^2(\nu_j d_{jk})} \le C_{\xi,\eta} \frac{(1 + \nu_j^2 d_{jk}^2)^{\nu/2}}{(1 + n d_{jk}^2)^{\nu/2}} = O\left(\left[\frac{\nu_j}{\sqrt{n}}\right]^{\nu}\right). \tag{5.36}$$

Hence, following (5.22) and denoting $u = \upsilon$ for $\gamma > 0$ and $u = \nu$ for $\gamma = 0$, we obtain

$$\begin{aligned} R_{121} &= O\left(\sum_{j=L}^{j_0} n^{-(u-1-\alpha_1)} 2^{2j[m_1(u-1-\alpha_1)-(r+1/2-1/p)]}\right) \\ &= O(n^{-2r/(2r+1)}) \end{aligned} \tag{5.37}$$

provided $\alpha_1$ satisfies (3.32) if $\gamma > 0$ (we have no restrictions on $\alpha$ if $\gamma > 0$). The last term, $R_{122}$, we partition into $R_{122} = R_{1221} + R_{1222}$ depending on the values of $\sqrt{n}|\theta_{jk}|$. Here

$$
\begin{aligned}
R_{1221} &= O\left(\sum_{j=L}^{j_0}\sum_{k=0}^{2^j-1} E[d_{jk}^2 \min(1, A_{jn}^2(d_{jk}))I(\sqrt{n}|\theta_{jk}| \le (\ln n)^b)]\right)\\
&= O\left(\sum_{j=L}^{j_0}\sum_{k=0}^{2^j-1} [n^{-1} + \theta_{jk}^2 I(\sqrt{n}|\theta_{jk}| \le (\ln n)^b)]\right)\\
&= O([n(\ln n)^{-2b}]^{-2r/(2r+1)})
\end{aligned}
\tag{5.38}
$$

by Lemma 7. For $R_{1222}$, repeating (5.24) and (5.25), we derive

$$
R_{1222} = O\left(\sum_{j=L}^{j_0}\sum_{k=0}^{2^j-1} n^{-1}[\eta_j^{-1}([\nu_j/\sqrt{n}]^{1+\alpha_1})I(\sqrt{n}|\theta_{jk}| > (\ln n)^b)]^2\right). \tag{5.39}
$$

By assumption on the decay of $\eta_j$ we derive that

$$
\eta_j^{-1}([\nu_j/\sqrt{n}]^{1+\alpha_1}) \le \begin{cases} C[\ln(n/\nu_j^2)]^{1/\gamma}, & \text{if } \gamma > 0,\\ C(n/\nu_j^2)^{(1+\alpha_1)/(2\nu)}, & \text{if } \gamma = 0. \end{cases} \tag{5.40}
$$

Consider the two cases separately. If $\gamma > 0$, then by (5.40) and Lemma 7,

$$
\begin{aligned}
R_{1222} &= O\left(\frac{(\ln n)^{2/\gamma}}{n}\sum_{j=L}^{j_0}\sum_{k=0}^{2^j-1} I(\sqrt{n}|\theta_{jk}| > (\ln n)^b)\right)\\
&= O(n^{-2r/(2r+1)}[\ln n]^{2/\gamma - pb}).
\end{aligned}
\tag{5.41}
$$

Now, to obtain (3.34), choose $b$ minimizing $\max(4br/(2r+1), 2/\gamma - pb)$. If $\gamma = 0$, then using (3.32) and (5.40), we derive by direct calculation that

$$
\begin{aligned}
R_{1222} &= O\left(\sum_{j=L}^{j_0} n^{-1}2^j(n/\nu_j^2)^{(1+\alpha_1)/\nu}\right)\\
&= O(n^{-2r/(2r+1)+(1+\alpha_1)/(\nu(2r+1))(1/p-1/2)}),
\end{aligned}
\tag{5.42}
$$

which agrees with (3.34).

*High resolution levels*: $j_1 + 1 \le j \le J-1$. Repeat (5.27) with $j_0$ replaced by $j_1$. Then $R_2 = O(n^{-2r/(2r+1)})$ follows from (3.29) and (5.27) and from the fact that $\sum_{j=j_1+1}^{J-1}\sum_{k=0}^{2^j-1} n/\nu_j^4 = O(n^{-2r/(2r+1)})$.

*Medium resolution levels*: $j_0 < j \le j_1$. Partition $R_5 = R_6 + R_7$ where

$$
R_6 = \sum_{j=j_0+1}^{j_1}\sum_{k=0}^{2^j-1} E(\hat{\theta}_{jk} - \theta_{jk})^2 I(\sqrt{n}|\theta_{jk}| > (\ln n)^b), \tag{5.43}
$$

and $R_7$ has $\leq$ instead of $>$ inside the indicator function; then further partition $R_6$ into $R_{61}$ and $R_{62}$ in a manner similar to (5.17). Repeat (5.18) and note that $m_2$ is chosen so that

$$\sum_{j=j_0+1}^{j_1} \sum_{k=0}^{2^j-1} [\nu_j^2 n^{-2} + \theta_{jk}^2 \nu_j^4 n^{-2}] = O(n^{-2r/(2r+1)}). \tag{5.44}$$

Then, by Lemma 7,

$$\begin{aligned} R_{61} &= O\Biggl(\sum_{j=j_0+1}^{j_1} \sum_{k=0}^{2^j-1} \biggl[\frac{\nu_j^2}{n^2} + \frac{\theta_{jk}^2 \nu_j^4}{n^2} + \frac{\sigma^2 I(\sqrt{n}|\theta_{jk}| > (\ln n)^b)}{n}\biggr]\Biggr) \\ &= O(n^{-2r/(2r+1)}). \end{aligned} \tag{5.45}$$

Partitioning $R_{62}$ further into $R_{621}$ and $R_{622}$ in a manner similar to $R_{121}$ and $R_{122}$ in (5.20), we derive

$$\begin{aligned} R_{621} &= O\Biggl(\sum_{j=j_0+1}^{j_1} \sum_{k=0}^{2^j-1} E\biggl[\frac{I_1(d_{jk})}{I_0(d_{jk})} \frac{A_{jn}(d_{jk}) I(\nu_j |d_{jk}| \leq C_\delta)}{1 + A_{jn}(d_{jk})}\biggr]^2 \\ &\qquad\qquad \times I(\sqrt{n}|\theta_{jk}| > (\ln n)^b)\Biggr) \\ &= O(n^{-2r/(2r+1)}) + O\Biggl(\sum_{j=j_0+1}^{j_1} \beta_{jn}^2 \nu_j^{2+2\delta} n^{-(1+\delta)} 2^{-2j(r+1/2-1/p)}\Biggr) \\ &= O(n^{-2r/(2r+1)}) \end{aligned} \tag{5.46}$$

by Lemmas 1 and 7, (5.21), condition (3.32) on $\alpha_2$ and the choice of $m_2$. For $R_{622}$ we obtain an upper bound similar to (5.39):

$$R_{622} = O\Biggl(\sum_{j=j_0+1}^{j_1} \sum_{k=0}^{2^j-1} n^{-1} [\eta_j^{-1}([\nu_j/\sqrt{n}]^{1+\alpha})]^2 I(\sqrt{n}|\theta_{jk}| > (\ln n)^b)\Biggr). \tag{5.47}$$

In the case where $\eta_j(x)$ has exponential decay ($\gamma > 0$) we just repeat (5.41) to obtain (3.34). In the case $\gamma = 0$ denote $(1 + \alpha_2)/\nu = h$ and observe that

$$\begin{aligned} R_{622} &= O\Biggl(\sum_{j=j_0+1}^{j_1} \sum_{k=0}^{2^j-1} n^{-1} (n/\nu_j^2)^h I(\sqrt{n}|\theta_{jk}| > (\ln n)^b)\Biggr) \\ &= O\Biggl(\sum_{j=j_0+1}^{j_1} \frac{n^{p/2-1+h}}{(\ln n)^{pb}} \frac{\nu_j^{-2h} \|\theta_{j\cdot}\|_{rp\infty}^p}{2^{2j(r+1/2-1/p)}}\Biggr) \\ &= O(n^{-2r/(2r+1)+h(1/p-1/2)(r+1)/(r(r+1/2))}), \end{aligned} \tag{5.48}$$

which agrees with (3.34). Now, to complete the proof we need to consider the term $R_7$ given by (5.43). Since

$$
\begin{aligned}
|\hat\theta_{jk}-\theta_{jk}| &\le \frac{|I_1(d_{jk})/I_0(d_{jk})-d_{jk}|}{1+A_{jn}(d_{jk})}+\frac{|d_{jk}-\theta_{jk}|}{1+A_{jn}(d_{jk})}\\
&\quad+\frac{|\theta_{jk}|}{1+A_{jn}(d_{jk})}+|\theta_{jk}|\\
&=O\biggl(\frac{\nu_j}{n}+\frac{\nu_j^2|d_{jk}|n^{-1}}{1+A_{jn}(d_{jk})}+\frac{|d_{jk}-\theta_{jk}|}{1+A_{jn}(d_{jk})}+|\theta_{jk}|\biggr)\\
&=O\biggl(|\theta_{jk}|+\frac{\nu_j}{n}+\frac{|d_{jk}-\theta_{jk}|}{1+A_{jn}(d_{jk})}\biggr),
\end{aligned}\tag{5.49}
$$

we partition $R_7$ into $R_{71}$, $R_{72}$ and $R_{73}$. Here, by Lemma 7 and (5.44)

$$
\begin{aligned}
R_{71}&=O\biggl(\sum_{j=j_0+1}^{j_1}\sum_{k=0}^{2^j-1}\theta_{jk}^2 I(\sqrt{n}|\theta_{jk}|\le(\ln n)^b)\biggr)\\
&=O([n(\ln n)^{-2b}]^{-2r/(2r+1)}),
\end{aligned}\tag{5.50}
$$

$$
R_{72}=O\biggl(\sum_{j=j_0+1}^{j_1}\sum_{k=0}^{2^j-1}\nu_j^2/n^2\biggr)=(n^{-2r/(2r+1)}).\tag{5.51}
$$

The third term, $R_{73}$, we partition into $R_{731}$ and $R_{732}$ where

$$
\begin{aligned}
R_{731}&=O\biggl(\sum_{j=j_0+1}^{j_1}\sum_{k=0}^{2^j-1}E\biggl[\frac{(d_{jk}-\theta_{jk})^2}{(1+A_{jn}(d_{jk}))^2}\biggr]I(1\le\sqrt{n}|\theta_{jk}|\le(\ln n)^b)\biggr),\\
R_{732}&=O\biggl(\sum_{j=j_0+1}^{j_1}\sum_{k=0}^{2^j-1}E\biggl[\frac{(d_{jk}-\theta_{jk})^2}{(1+A_{jn}(d_{jk}))^2}\biggr]I(\sqrt{n}|\theta_{jk}|\le 1)\biggr).
\end{aligned}\tag{5.52}
$$

By Lemma 7,

$$
\begin{aligned}
R_{731}&=O\biggl(\sum_{j=j_0+1}^{j_1}\sum_{k=0}^{2^j-1}n^{-1}I(n^{-1}<\theta_{jk}^2)I(\sqrt{n}|\theta_{jk}|\le(\ln n)^b)\biggr)\\
&=O(n^{-2r/(2r+1)}[\ln n]^{4br/(2r+1)}).
\end{aligned}\tag{5.53}
$$

For an upper bound for $R_{732}$ note that since $\xi(x)$ is bounded

$$
\begin{aligned}
&E[(d_{jk}-\theta_{jk})^2(1+A_{jn}(d_{jk}))^{-2}]I(\sqrt{n}|\theta_{jk}|\le 1)\\
&\quad=O\biggl(\int_{-\infty}^{\infty}\frac{(x-\theta_{jk})^2}{(1+\beta_{jn}\sqrt{n}\nu_j^{-1}\eta_j(\sqrt{n}x))^2}\\
&\qquad\qquad\times\frac{\sqrt{n}}{\sqrt{2\pi}\sigma}e^{-n(x-\theta_{jk})^2/2\sigma^2}\,dx\,I(\sqrt{n}|\theta_{jk}|\le 1)\biggr)
\end{aligned}
$$

$$= O\left(\nu_j^2 \beta_{jn}^{-2} n^{-2} \int_{-\infty}^{\infty} z^2 \exp(-z^2/2\sigma^2)[\eta_j(1+|z|)]^{-2}\,dz\right)$$
$$= O(\nu_j^2 \beta_{jn}^{-2} n^{-2})$$

by condition (A5). Hence, whenever $\alpha_2 \geq 0$, by (5.44) we have

$$R_{732} = O\left(\sum_{j=j_0+1}^{j_1} \sum_{k=0}^{2^j-1} \nu_j^2 n^{-2}\right) = O(n^{-2r/(2r+1)}). \tag{5.54}$$

Combination of (5.37), (5.38), (5.41), (5.42) and (5.45)–(5.54) completes the proof. $\square$

Proof of Theorem 6. Note that for $j \leq j_1$, the derivation of the error is identical to that for Theorem 5. When $j > j_1$, the proof is very similar to the proof of Theorem 3 with $j_0$ replaced by $j_1$ and $\lambda_\eta = 1$. Again, since $|\hat{\theta}_{jk} - \theta_{jk}| \leq |\theta_{jk}| + |\hat{\theta}_{jk}|$ and due to (3.29), we shall be concerned with the $|\hat{\theta}_{jk}|$ term only. Noting that $|\hat{\theta}_{jk}| \leq |I_1(d_{jk})/I_0(d_{jk})|$ and using the inequality (5.28), we partition the error and derive upper bounds identical to those in (5.29), (5.30) and (5.32), which leads to $\sum_{j_1+1}^{J_0} \sum_{k=0}^{2^j-1} E(\hat{\theta}_{jk} - \theta_{jk})^2 = O(n^{-2r/(2r+1)}(\ln n)^{1/(2r+1)})$. Derivation of the upper bound for $\sum_{j=J_0+1}^{J-1} \sum_{k=0}^{2^j-1} E(\hat{\theta}_{jk} - \theta_{jk})^2$ repeats the derivation of $R_4$ in Theorem 3, and similarly, we obtain $\sum_{j=J_0+1}^{J-1} \sum_{k=0}^{2^j-1} E(\hat{\theta}_{jk} - \theta_{jk})^2 = O(n^{-2r/(2r+1)})$. $\square$

Proof of Theorem 7. The validity of Theorem 6 follows from Theorems 1, 3, 5 and 6. The only portion of the error which needs additional consideration is the one corresponding to $j > J_0$. Recall that

$$\hat{\theta}_{jk} = \frac{\lambda_j^* I_1(d_{jk}) + (1-\lambda_j^*) I_1^*(d_{jk})}{\lambda_j^* I_0(d_{jk}) + (1-\lambda_j^*) I_0^*(d_{jk}) + \beta_{jn}\sqrt{n}\eta_j(\sqrt{n} d_{jk})}, \tag{5.55}$$

where $\eta_j(x)$ and $\lambda_j^*$ are given by (3.36),

$$I_i(d) = \int_{-\infty}^{\infty} x^i \sqrt{n} (\sqrt{2\pi}\sigma_0)^{-1} \exp(-n(d-x)^2/(2\sigma_0^2)) \nu_j \xi(\nu_j x)\,dx, \qquad i = 0, 1,$$

$$I_i^*(d) = \int_{-\infty}^{\infty} x^i \sqrt{n} \zeta(\sqrt{n}(d-x)) \nu_j \xi(\nu_j x)\,dx, \qquad i = 0, 1.$$

Consider the cases when $\xi$ is a normal and a heavy-tailed prior separately. If $\xi$ is a normal prior, then by Lemma 2

$$\begin{aligned}
R_8 &= \sum_{j=J_0+1}^{J-1} \sum_{k=0}^{2^j-1} E(\hat{\theta}_{jk} - \theta_{jk})^2 \\
&= O\left( \sum_{j=J_0+1}^{J-1} \sum_{k=0}^{2^j-1} E\left[ \left( \frac{I_1(d_{jk})}{I_0(d_{jk})} \right)^2 + \left( \frac{I_1^*(d_{jk})}{I_0^*(d_{jk})} \right)^2 + \theta_{jk}^2 \right] \right) \qquad (5.56) \\
&= O\left( \sum_{j=J_0+1}^{J-1} \sum_{k=0}^{2^j-1} (n^2 \nu_j^{-4} E d_{jk}^2 + n\nu_j^{-4} + \theta_{jk}^2) \right) = O(n^{-2r/(2r+1)}).
\end{aligned}$$

If $\lambda_\xi = 0$, then partition $R_8$ as $R_8 = R_{81} + R_{82}$ where

$$R_{81} = O\left( \sum_{j=J_0+1}^{J-1} \sum_{k=0}^{2^j-1} \left[ E\left\{ \left( \frac{I_1(d_{jk})}{I_0(d_{jk})} \right)^2 I\left( \frac{n|d_{jk}|}{\nu_j} \to 0 \right) \right\} + E\left( \frac{I_1^*(d_{jk})}{I_0^*(d_{jk})} \right)^2 + \theta_{jk}^2 \right] \right),$$

so that $R_{81} = O(n^{-2r/(2r+1)})$ by Lemma 2 and calculation similar to (5.56), and

$$R_{82} = O\left( \sum_{j=J_0+1}^{J-1} \sum_{k=0}^{2^j-1} E\left[ \left( \frac{I_1(d_{jk})}{I_0(d_{jk}) + (\lambda_j^*)^{-1} \beta_{jn} \sqrt{n} \zeta(\sqrt{n} d_{jk})} \right)^2 \times I\left( \frac{n|d_{jk}|}{\nu_j} \geq M \right) \right] \right).$$

Let $\zeta_1(x)$ be a heavy-tailed p.d.f. such that

$$C_{\zeta,1} (\sigma_0 \sqrt{2\pi})^{-1} \exp\{-x^2/2\sigma_0^2\} \leq \zeta_1(x) \leq C_{\zeta,2} (x^2+1)^{-1} \zeta(x) \qquad (5.57)$$

for some positive $C_{\zeta,1}$ and $C_{\zeta,2}$. Then, it is easy to see that $I_0(d) \leq C_{\zeta,1}^{-1} \int \sqrt{n} \zeta_1 \times (\sqrt{n}(d-x)) \nu_j \xi(\nu_j x)\, dx \sim \sqrt{n} \zeta_1(\sqrt{n} d)$. Thus, by Lemma 2 $\sqrt{n}\zeta(\sqrt{n} d_{jk})/ I_0(d_{jk}) \geq C_{\zeta,1}^{-1} \zeta(\sqrt{n} d_{jk})/\zeta_1(\sqrt{n} d_{jk})$ and

$$\begin{aligned}
R_{82} &= O\left( \sum_{j=J_0+1}^{J-1} \sum_{k=0}^{2^j-1} E \frac{[I_1(d_{jk})/I_0(d_{jk})]^2 I(\nu_j^{-1} n |d_{jk}| \geq M)}{[\beta_{jn} (\lambda_j^*)^{-1} \zeta(\sqrt{n} d_{jk})/\zeta_1(\sqrt{n} d_{jk})]^2} \right) \\
&= O\left( \sum_{j=J_0+1}^{J-1} \sum_{k=0}^{2^j-1} E[(n d_{jk}^2)^{-1} d_{jk}^2 I(\nu_j^{-1} n |d_{jk}| \geq M)] \right) \qquad (5.58) \\
&= O(n^{-2r/(2r+1)}),
\end{aligned}$$

which completes the proof. □

Proof of Theorem 8. The proof is based on the following obvious corollary of Kolmogorov's three series theorem: if $Z_1, Z_2, \ldots, Z_n$ are independent variables such that $\sum_n E|Z_n| < \infty$, then the series $\sum_n Z_n$ converges with probability 1 (see, e.g., [22], Section 4.2). Since we impose priors on each coefficient independently, it follows from the above that $f \in B^r_{pq}$ if

$$S_1 = \sum_{j=L}^{\infty} 2^{j(r+1/2-1/p)q} E\left[\sum_{k=0}^{2^j-1} |\theta_{jk}|^p\right]^{q/p} < \infty \tag{5.59}$$

where the expectation is taken over the prior p.d.f. $\pi_{jn}\nu_j\xi(\nu_j x) + (1 - \pi_{jn})\delta(0)$. Note that for $a > 0$ we have $E|\theta_{jk}|^a = C_a \pi_{jn}\nu_j^{-a}$ where $C_a = \int |x|^a \xi(x)\,dx$ provided the last integral is convergent. If $q/p \le 1$, then using Jensen's inequality $E(|Y|^{q/p}) \le (E|Y|)^{q/p}$ and the expression for $E|\theta_{jk}|^q$ we obtain

$$S_1 \le \sum_{j=L}^{\infty} 2^{j(r+1/2-1/p)q} \left[\sum_{k=0}^{2^j-1} E|\theta_{jk}|^p\right]^{q/p} \le \sum_{j=L}^{\infty} 2^{j(r+1/2)q} \nu_j^{-q} \pi_{jn}^{q/p}. \tag{5.60}$$

Taking into account that $\beta_{jn} = (1 + \pi_{jn})^{-1}$, we derive (3.41).

If $q > p$, note that if $f \in B^r_{p,p}$, then $f \in B^r_{p,q}$, so it is enough to consider $q = p$. Repeating the previous calculation, we arrive at (3.41) again. □

Proof of Corollary 4. Plug the values of $\nu_j$ and $\beta_{jn}$ into (3.41) and check that the sum (3.41) is uniformly bounded by a constant independent of $n$. □

Proof of Corollary 5. Plugging the values of $\nu_j$ given by (3.28) into (3.41), we derive that (3.41) holds whenever the sum $S^* = \sum_{j=L}^{\infty} 2^{jv} \times \beta_{jn}^{-\min(q/p,1)}$ is uniformly bounded, where

$$v = \begin{cases} (1/p - 1/2)/2, & L \le j \le j_0, \\ (1/p - 1/2)(1 + 1/r), & j_0 < j \le j_1, \\ 0, & j > j_1. \end{cases} \tag{5.61}$$

Now to complete the proof, check that $S^*$ is uniformly bounded whenever (3.43) is valid. □

The proof of Theorem 9 is based on the following lemma.

Lemma 8. *Under the conditions of Theorem 9, for some positive absolute constant $C$*

$$R_n(B^r_{p,q}(A), \hat{f}\,) \ge C \sum_{j=L}^{j_0} 2^{-2rj} I\left(\beta_{jn} > \frac{\nu_j}{\sqrt{n}} \eta_j\left(\frac{C_\delta \sqrt{n}}{\nu_j}\right)\right) - O(n^{-2r/(2r+1)}). \tag{5.62}$$

PROOF. Repeating the beginning of the proof of Theorem 1, we derive that $R_n(B^r_{p,q}(A), \hat{f}) \geq R \geq R_1$ where $R_1$ is defined in (5.15). Since for any random variables $X$ and $Y$ $E(X+Y)^2 \geq 0.5EX^2 - EY^2$, we obtain that $R_1 \geq 0.5R_{12} - R_{11}$ where, just as in (5.19), $R_{11} = O(n^{-2r/(2r+1)})$. Now, by (5.20) $R_{12} \geq R_{121}$. Consequently,

$$R_n(B^r_{p,q}(A), \hat{f}\,) \geq 0.5R_{121} - O(n^{-2r/(2r+1)}) \tag{5.63}$$

where $R_{121}$ is defined in (5.20). Recall (5.23) and by direct calculation verify that $\alpha_{jn} \geq C_\delta \sqrt{n}\nu_j^{-1}$ if and only if

$$\beta_{jn} \geq \nu_j[\sqrt{n}\eta(C_\delta\sqrt{n}\nu_j^{-1})]^{-1}. \tag{5.64}$$

For the values of $j$ for which (5.64) holds, $A^2_{jn}(d_{jk}) \geq 1$, so that $[A_{jn}(d_{jk})/(1+A_{jn}(d_{jk}))]^2 \geq 0.25$. Hence

$$\begin{aligned} R_{121} &\geq 0.25 \sum_{j=L}^{j_0} \sum_{k=0}^{2^j-1} E[d^2_{jk} I(C_\delta \nu_j^{-1} < |d_{jk}| \leq \alpha_{jn}/\sqrt{n}\,)] \\ &\geq 0.25 C_\delta^2 \sum_{j=L}^{j_0} 2^j \nu_j^{-2} I(\beta_{jn} \geq \nu_j[\sqrt{n}\eta(C_\delta\sqrt{n}\nu_j^{-1})]^{-1}). \end{aligned} \tag{5.65}$$

Plug the value of $\nu_j$ into (5.65). The combination of (5.63)–(5.65) completes the proof. Note that in Theorem 1 the choice of $\beta_{jn}$ makes the inequality $\alpha_{jn} \geq C_\delta\sqrt{n}\nu_j^{-1}$ impossible. □

PROOF OF THEOREM 9. Plugging $\beta_{jn} = 2^{bj}$ into (5.62) we obtain that

$$\beta_{jn} > \nu_j[\sqrt{n}\eta(C_\delta\sqrt{n}/\nu_j)]^{-1} \quad \Longleftrightarrow \quad \nu_j \geq C_\delta\sqrt{n}[\eta^{-1}(2^{j(r+1/2-b)}n^{-1/2})]^{-1}.$$

Recalling that $j \leq j_0$, we note that $2^{j(r+1/2-b)}n^{-1/2} \leq n^{-b/(2r+1)}$. Hence, by Lemma 5.62,

$$\begin{aligned} &R_n(B^r_{p,q}(A), \hat{f}\,) \\ &\quad \geq C \sum_{j=L}^{j_0} 2^{-2rj} I(\beta_{jn} > \nu_j[\sqrt{n}\eta(C_\delta\sqrt{n}/\nu_j)]^{-1}) \\ &\quad \geq C \sum_{j=L}^{j_0} 2^{-2rj} I(2^j \geq [C_\delta\sqrt{n}\,]^{1/(r+1/2)}[\eta^{-1}(n^{-b/(2r+1)})]^{-1/(r+1/2)}) \\ &\quad = C_1 n^{-2r/(2r+1)}[\eta^{-1}(n^{-b/(2r+1)})]^{2r/(r+1/2)} \end{aligned}$$

for some positive absolute constants $C$ and $C_1$, which proves (3.45). To finish the proof, consider various cases for $\eta(x)$. □

**Acknowledgment.** The author thanks the referees for their extremely helpful and constructive comments.

Department of Mathematics
University of Central Florida
Orlando, Florida 32816-1364
USA
E-mail: mpensky@pegasus.cc.ucf.edu